\documentclass[english]{svjour3}
\usepackage[T1]{fontenc}
\usepackage[latin9]{inputenc}
\usepackage{units}
\usepackage{url}
\usepackage{amsmath}
\usepackage{amssymb}
\usepackage{graphicx}
\usepackage{esint}

\makeatletter

\providecommand{\tabularnewline}{\\}

  \newenvironment{svmultproof}{\begin{proof}}{\qed\end{proof}}

\RequirePackage{fix-cm}

\smartqed  

\usepackage[numbers]{natbib}

\makeatother

\usepackage{babel}

\begin{document}

\title{Bounds on Tail Probabilities in Exponential families}

\titlerunning{Bounds on Tail Probabilities}

\author{Peter Harremo{\"e}s}

\institute{Peter Harremo{\"e}s\at Copenhagen Business College\\
Copenhagen, Denmark\\
Tel.: +45 39 56 41 71\\
\email{harremoes@ieee.org}\\
}

\date{Received: date / Accepted: date}
\maketitle
\begin{abstract}
In this paper we present various new inequalities for tail proabilities
for distributions that are elements of the most improtant exponential
families. These families include the Poisson distributions, the Gamma
distributions, the binomial distributions, the negative binomial distributions
and the inverse Gaussian distributions. All these exponential families
have simple variance functions and the variance functions play an
important role in the exposition. All the inequalities presented in
this paper are formulated in terms of the signed log-likelihood. The
inequalities are of a qualitative nature in that they can be formulated
either in terms of stochastic domination or in terms of an intersection
property that states that a certain discrete distribution is very
close to a certain continuous distribution.

\end{abstract}

\section{Introduction}

Let $X_{1},\dots,X_{n}$ be i.i.d. random variables such that the
moment generating function $\beta\curvearrowright E[\exp(\beta X_{1})]$
is finite in a neighborhood of the origin. For a fixed value of $\mu$
one is interested in approximating the tail distribution: $\Pr(\sum_{i=1}^{n}X_{i}\leq\mu)$
. If $\mu$ is close to the mean of $X_{1}$ one would usually approximate
the tail probability by the tail probability of a Gaussian random
variable. If $\mu$ is far from the mean of $X_{1}$ the tail probability
can be estimated using large deviation theory. According to the Sanov
theorem the probability that the deviation from the mean is as large
as $\mu$ is of the order $\exp\left(-D\right)$ where D is a constant
or to be more precise
\[
-\frac{\ln\left(\Pr(\sum_{i=1}^{n}X_{i}\leq\mu)\right)}{n}\to D
\]
for $n\to\infty.$ Bahadur and Rao \cite{Bahadur1960,Bahadur1960a}
improved the estimate of this large deviation probability, and in
\cite{Gyorfi2012} such Gaussian tail approximations were extended
to situations where one normally uses large deviation techniques.

The distribution of the signed log-likelihood is close to a standard
Gaussian for a variety of distributions. An asymptotic result for
large sample sizes this is not new \cite{Zhang2008}, but in this
paper we are interested in inequalities that hold for any sample size.
Some inequalities of this type can be found in \cite{Alfers1984,Reiczigel2011,Harremoes2012,Zubkov2013,Harremoes2014},
but here we attempt to give a more systematic presentation including
a number of new or improved inequalities. 

In this paper we let $\tau$ denote the circle constant $2\pi$ and
$\phi$ will denote the standard Gaussian density
\[
\frac{\exp\left(-\frac{x^{2}}{2}\right)}{\tau^{\nicefrac{1}{2}}}.
\]
 We let $\Phi$ denote the distribution function of the standard Gaussian
\[
\Phi\left(t\right)=\int_{-\infty}^{t}\phi\left(x\right)\,\text{d}x~.
\]

The rest of the paper is organized as follows. In Section \ref{sec:signedExpo}
we define the signed log-likelihood of exponential families and look
at some of the fundamental properties of the signed log-likelihood.
Next we study inequalities for the signed log-likelihood for certain
exponential families associated with continuous waiting times. We
start with the inverse Gaussian in Section \ref{sec:Inequalities-for-Wald}
that is particularly simple. Then we study the exponential distributions
(Section \ref{sec:Exponential-distributions}) and more general Gamma
distributions (Section \ref{sec:Gamma-distributions}). Next we turn
our attention to discrete waiting times. First we obtain some new
inequalities for the geometric distributions (Section \ref{sec:Geometric-distributions})
and then we generalize the results to negative binomial distributions
(Section \ref{sec:Inequalities-for-waiting}). The negative binomial
distributions are waiting times in Bernoulli processes, so in Section
\ref{sec:Inequalities-for-binomial} our inequalities between negative
binomial distributions and Gamma distributions are translated into
inequalities between binomial distributions and Poisson distributions.
Combined with our domination inequalities for Gamma distributions
we obtain an intersection inequality between binomial distributions
and the Standard Gaussian distribution. In this paper the focus is
on intersection inequalities and stochastic domination inequalities,
but in the discussion we mention some related inequalities of other
types and how they may be improved.

\section{The signed log-likelihood for exponential families\label{sec:signedExpo}}

Consider the 1-dimensional exponential family $P_{\beta}$ where
\[
\frac{\text{d}P_{\beta}}{\text{d}P_{0}}\left(x\right)=\frac{\exp\left(\beta\cdot x\right)}{Z\left(\beta\right)}
\]
and $Z$ denotes the \emph{moment generating function} given by Z$\left(\beta\right)=\int\exp\left(\beta\cdot x\right)\,\text{d}P_{0}x.$
Let $P^{\mu}$ denote the element in the exponential family with mean
value $\mu,$ and let $\hat{\beta}\left(\mu\right)$ denote the corresponding
maximum likelihood estimate of $\beta.$ Let $\mu_{0}$ denote the
mean value of $P_{0}$. Then
\begin{align*}
D\left(P^{\mu}\Vert P_{0}\right) & =\int\ln\left(\frac{\text{d}P^{\mu}}{\text{d}P_{0}}\left(x\right)\right)\,\text{d}P^{\mu}x.
\end{align*}
The variance function of an exponential family is defined so that
$V\left(\mu\right)$ is the variance of $P^{\mu}.$ The variance functions
uniquely characterizes the corresponding exponential families and
most important exponential families have very simple variance functions.
If we know the varince function the divergence can be calculated according
to the following formula.
\begin{lemma}
In an exponential family $\left(P^{\mu}\right)$ parametrized by mean
value $\mu$and with variance function $V$ information divergence
can be calculated according to the formula
\[
D\left(P^{\mu_{1}}\Vert P^{\mu_{2}}\right)=\int_{\mu_{1}}^{\mu_{2}}\frac{\mu-\mu_{1}}{V\left(\mu\right)}\,\mathrm{d}\mu.
\]
\end{lemma}
\begin{svmultproof}
The divergence is given by 
\begin{align*}
D\left(P^{\mu_{1}}\Vert P^{\mu_{2}}\right) & =\int\ln\left(\frac{\text{d}P^{\mu_{1}}}{\text{d}P^{\mu_{2}}}\left(x\right)\right)\,\text{d}P^{\mu_{1}}x\\
 & =\int\ln\left(\frac{\frac{\exp\left(\beta_{1}\cdot x\right)}{Z\left(\beta_{1}\right)}}{\frac{\exp\left(\beta_{2}\cdot x\right)}{Z\left(\beta_{2}\right)}}\right)\,\text{d}P^{\mu_{1}}x\\
 & =\int\left(\left(\beta_{1}-\beta_{2}\right)x-\ln Z\left(\beta_{1}\right)+\ln Z\left(\beta_{2}\right)\right)\,\text{d}P^{\mu_{1}}x\\
 & =\left(\beta_{1}-\beta_{2}\right)\mu_{1}-\ln Z\left(\beta_{1}\right)+\ln Z\left(\beta_{2}\right).
\end{align*}
The derivative with respect to $\beta_{2}$ is 
\[
\frac{\partial}{\partial\beta_{2}}D\left(P^{\mu_{1}}\Vert P^{\mu_{2}}\right)=\mu_{2}-\mu_{1}.
\]
Therefore the derivative with respect to $\mu_{2}$ is
\begin{align*}
\frac{\partial}{\partial\mu_{2}}D\left(P^{\mu_{1}}\Vert P^{\mu_{2}}\right) & =\frac{\mu_{2}-\mu_{1}}{\frac{\mathrm{d}\mu_{2}}{\mathrm{d}\beta_{2}}}\\
 & =\frac{\mu_{2}-\mu_{1}}{V\left(\mu_{2}\right)}.
\end{align*}
Together with the trivial identity 
\[
D\left(P^{\mu_{1}}\Vert P^{\mu_{1}}\right)=\int_{\mu_{1}}^{\mu_{1}}\frac{\mu-\mu_{1}}{V\left(\mu\right)}\,\mathrm{d}\mu
\]
the results follows.\end{svmultproof}

\begin{definition}
(From \cite{Barndorff-Nielsen1990}) Let $X$ be a random variable
with distribution $P_{0}.$ Then the \emph{signed log-likelihood}
$G\left(X\right)$ of $X$ is the random variable given by
\[
G\left(x\right)=\left\{ \begin{array}{cc}
-\left[2D\left(P^{x}\Vert P_{0}\right)\right]^{\nicefrac{1}{2}}, & \text{for }x<\mu_{0};\\
+\left[2D\left(P^{x}\Vert P_{0}\right)\right]^{\nicefrac{1}{2}}, & \text{for }x\geq\mu_{0}.
\end{array}\right.
\]

\end{definition}
We will need the following general lemma.
\begin{lemma}
\label{lem:Voksende}If the variance function is increasing then
\[
\frac{G\left(x\right)}{x-\mu_{0}}
\]
is a decreasing function of $x$.\end{lemma}
\begin{svmultproof}
We have 
\begin{align*}
\frac{\mathrm{d}}{\mathrm{d}x}\left(\frac{G\left(x\right)}{x-\mu_{0}}\right) & =\frac{\left(x-\mu_{0}\right)\frac{D'\left(x\right)}{G\left(x\right)}-G\left(x\right)}{\left(x-\mu_{0}\right)^{2}}\\
 & =\frac{\left(x-\mu_{0}\right)\int_{x}^{\mu_{0}}\frac{-1}{V\left(\mu\right)}\,\mathrm{d}\mu-2D}{\left(x-\mu_{0}\right)^{2}G\left(x\right)}\\
 & =\frac{\left(x-\mu_{0}\right)\int_{\mu_{0}}^{x}\frac{1}{V\left(\mu\right)}\,\mathrm{d}\mu-2D}{\left(x-\mu_{0}\right)^{2}G\left(x\right)}.
\end{align*}
We have to prove that numerator is positive for $x<\mu_{0}$ and negative
for $x>\mu_{0}.$ The numerator can be calculated as 
\begin{align*}
\left(x-\mu_{0}\right)\int_{\mu_{0}}^{x}\frac{1}{V\left(\mu\right)}\,\mathrm{d}\mu-2D & =\left(x-\mu_{0}\right)\int_{\mu_{0}}^{x}\frac{1}{V\left(\mu\right)}\,\mathrm{d}\mu-2\int_{\mu_{0}}^{x}\frac{\mu-\mu_{0}}{V\left(\mu\right)}\,\mathrm{d}\mu\\
 & =\int_{\mu_{0}}^{x}\left(\frac{x-\mu_{0}}{V\left(\mu\right)}-2\frac{\mu-\mu_{0}}{V\left(\mu\right)}\right)\,\mathrm{d}\mu\\
 & =\int_{\mu_{0}}^{x}\frac{x+\mu_{0}-2\mu}{V\left(\mu\right)}\,\mathrm{d}\mu
\end{align*}
 If $x>\mu_{0}$ then
\begin{align*}
\int_{\mu_{0}}^{x}\frac{x+\mu_{0}-2\mu}{V\left(\mu\right)}\,\mathrm{d}\mu & =\int_{\mu_{0}}^{\frac{x+\mu_{0}}{2}}\frac{x+\mu_{0}-2\mu}{V\left(\mu\right)}\,\mathrm{d}\mu+\int_{\frac{x+\mu_{0}}{2}_{0}}^{x}\frac{x+\mu_{0}-2\mu}{V\left(\mu\right)}\,\mathrm{d}\mu\\
 & \geq\int_{\mu_{0}}^{\frac{x+\mu_{0}}{2}}\frac{x+\mu_{0}-2\mu}{V\left(\frac{x+\mu_{0}}{2}\right)}\,\mathrm{d}\mu+\int_{\frac{x+\mu_{0}}{2}_{0}}^{x}\frac{x+\mu_{0}-2\mu}{V\left(\frac{x+\mu_{0}}{2}\right)}\,\mathrm{d}\mu\\
 & =\int_{\mu_{0}}^{x}\frac{x+\mu_{0}-2\mu}{V\left(\frac{x+\mu_{0}}{2}\right)}\,\mathrm{d}\mu=0.
\end{align*}
The inequality for $x<\mu_{0}$ is proved in the same way.
\end{svmultproof}

\section{Inequalities for inverse Gaussian\label{sec:Inequalities-for-Wald}}

The \emph{inverse Gaussian distribution} and it is used to model waiting
times for a Wiener process (Brownian motion) with drift. An inverse
Gaussian distribution has density 
\[
f\left(w\right)=\left[\frac{\lambda}{\tau w^{3}}\right]^{\nicefrac{1}{2}}\exp\left(\frac{-\lambda\left(w-\mu\right)^{2}}{2\mu^{2}w}\right)
\]
with mean value parameter $\mu$ and shape parameter $\lambda.$ The
variance function is $V\left(\mu\right)=\mu^{3}/\lambda.$ 

The divergence of an inverse Gaussian distribution with mean $\mu_{1}$
from an inverse Gaussian distribution with mean $\mu_{2}$ is
\[
\int_{\mu_{1}}^{\mu_{2}}\frac{\mu-\mu_{1}}{\mu^{3}/\lambda}\,d\mu=\frac{\lambda\left(\mu_{1}-\mu_{2}\right)^{2}}{2\mu_{1}\mu_{2}^{2}}.
\]
Hence the signed log-likelihood is 
\[
G_{\mu,\lambda}\left(w\right)=\frac{\lambda^{\nicefrac{1}{2}}\left(w-\mu\right)}{w^{\nicefrac{1}{2}}\mu}.
\]
We observe that 
\begin{align*}
G_{\mu,\lambda}\left(w\right) & =\left[\frac{\lambda}{\mu}\right]^{\nicefrac{1}{2}}\frac{\frac{w}{\mu}-1}{\left[\frac{w}{\mu}\right]^{\nicefrac{1}{2}}}\\
 & =\left[\frac{\lambda}{\mu}\right]^{\nicefrac{1}{2}}G_{1,1}\left(\frac{w}{\mu}\right).
\end{align*}
Note that the \emph{saddle-point approximation} \cite{Daniels1954}
is exact for the family of inverse Gaussian distributions, i.e.
\[
f\left(w\right)=\frac{\phi\left(G\left(w\right)\right)}{\left[V\left(w\right)\right]^{\nicefrac{1}{2}}}.
\]

\begin{figure}
\begin{centering}
\includegraphics[scale=0.8]{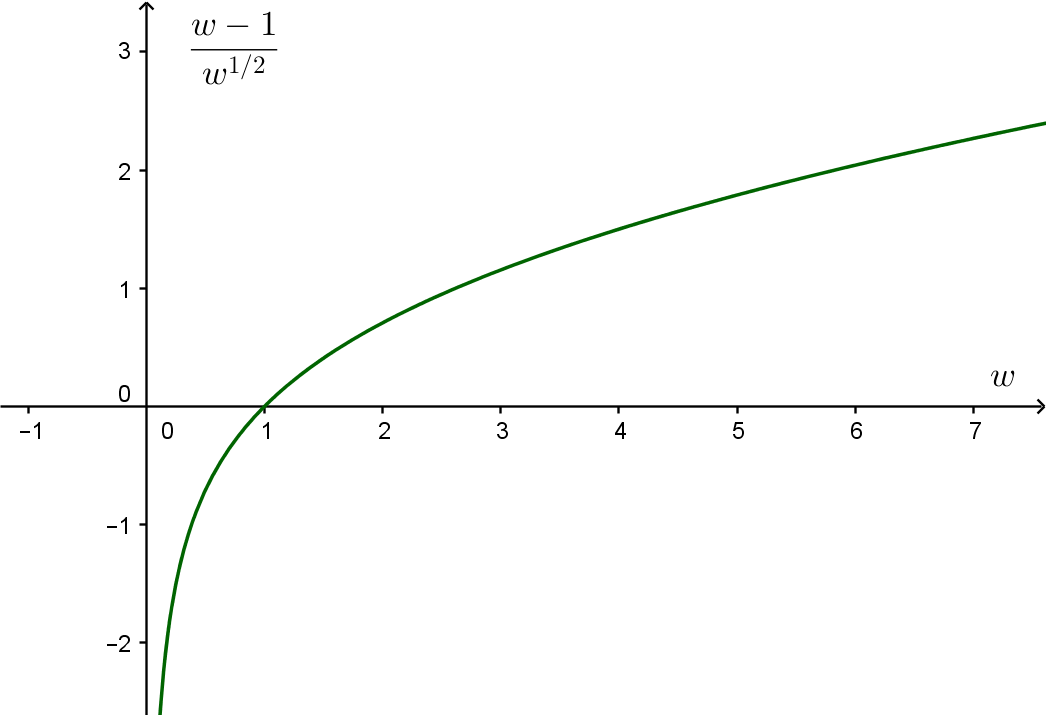}
\par\end{centering}

\caption{The signed log-likelihood of an inverse Gaussian distribution with
mean value 1 and shape parameter 1. }
\end{figure}

\begin{lemma}
The probability density of the random variable $G_{\mu,\lambda}\left(W\right)$
is 
\[
\frac{2\phi\left(z\right)}{1+g^{-1}\left(z\cdot\left[\frac{\mu}{\lambda}\right]^{\nicefrac{1}{2}}\right)}
\]
where $g$ denotes the function $G_{1,1}$.\end{lemma}
\begin{svmultproof}
The density of $G_{\mu,\lambda}\left(W\right)$ is 
\begin{align*}
\frac{f\left(G_{\mu,\lambda}^{-1}\left(z\right)\right)}{G_{\mu,\lambda}^{'}\left(G_{\mu,\lambda}^{-1}\left(z\right)\right)} & =\frac{\frac{\phi\left(G_{\mu,\lambda}\left(G_{\mu,\lambda}^{-1}\left(z\right)\right)\right)}{\left[V\left(G_{\mu,\lambda}^{-1}\left(z\right)\right)\right]^{\nicefrac{1}{2}}}}{G_{\mu,\lambda}^{'}\left(G_{\mu,\lambda}^{-1}\left(z\right)\right)}\\
 & =\frac{\phi\left(z\right)}{\left[V\left(G_{\mu,\lambda}^{-1}\left(z\right)\right)\right]^{\nicefrac{1}{2}}G_{\mu,\lambda}^{'}\left(G_{\mu,\lambda}^{-1}\left(z\right)\right)}.
\end{align*}
Now we use that 
\begin{align*}
G_{\mu,\lambda}^{'}\left(w\right) & =\frac{w^{\nicefrac{1}{2}}\mu\lambda^{\nicefrac{1}{2}}-\nicefrac{1}{2}\cdot w^{-\nicefrac{1}{2}}\mu\lambda^{\nicefrac{1}{2}}\left(w-\mu\right)}{w\mu^{2}}\\
 & =\lambda^{\nicefrac{1}{2}}\frac{\mu+w}{2w^{\nicefrac{3}{2}}\mu}.
\end{align*}
Hence
\begin{align*}
\left[V\left(w\right)\right]^{\nicefrac{1}{2}}G_{\mu,\lambda}^{'}\left(w\right) & =\left[\frac{w^{3}}{\lambda}\right]^{\nicefrac{1}{2}}\cdot\lambda^{\nicefrac{1}{2}}\frac{\mu+w}{2w^{\nicefrac{3}{2}}\mu}\\
 & =\frac{\mu+w}{2\mu}.
\end{align*}
Therefore the density of $G_{\mu,\lambda}\left(W\right)$ is 
\begin{align*}
f\left(z\right) & =\frac{\phi\left(z\right)2\mu}{G_{\mu,\lambda}^{-1}\left(z\right)+\mu}.
\end{align*}
By isolating $x$ in the equation $G_{\mu,\lambda}\left(x\right)=z$
we get 
\[
G_{\mu,\lambda}^{-1}\left(z\right)=\mu\cdot G_{1,1}^{-1}\left(z\cdot\left[\frac{\mu}{\lambda}\right]^{\nicefrac{1}{2}}\right).
\]
 Hence
\[
f\left(z\right)=\frac{\phi\left(z\right)2}{G_{1,1}^{-1}\left(z\cdot\left[\frac{\mu}{\lambda}\right]^{\nicefrac{1}{2}}\right)+1}\,,
\]
which proves the theorem.
\end{svmultproof}

\begin{figure}
\begin{centering}
\includegraphics[scale=1.2]{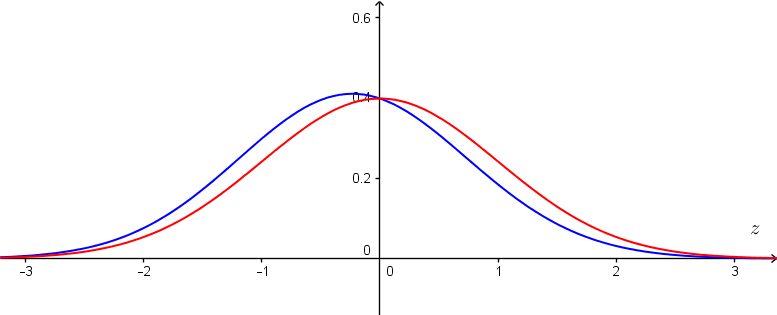}
\par\end{centering}

\caption{The density of the signed log-likelihood of an inverse Gaussian distribution
(blue) with mean value 1 and shape parameter 1 comparead with the
density of a standard Gaussian distribution (red). }
\end{figure}

\begin{lemma}
\label{lem:dominans} (From \cite{Harremoes2012}) Let $X_{1}$ and
$X_{2}$ denote random variables with density functions $f_{1}$ and
$f_{2}$. If $f_{1}\left(x\right)\geq f_{2}\left(x\right)$ for $x\leq x_{0}$
and $f_{1}\left(x\right)\leq f_{2}\left(x\right)$ for $x\geq x_{0},$
then $X_{1}$ is stochastically dominated by $X_{2}.$ In particular
if $\frac{f_{2}\left(x\right)}{f_{1}\left(x\right)}$ is increasing
then $X_{1}$ is stochastically dominated by $X_{2}.$ \end{lemma}
\begin{svmultproof}
Assume that $f_{1}\left(x\right)\geq f_{2}\left(x\right)$ for $x\leq x_{0}$
and that $f_{1}\left(x\right)\leq f_{2}\left(x\right)$ for $x\geq x_{0}.$
For $t\geq x_{0}$ we have 
\begin{align*}
P\left(X_{1}\geq t\right) & =\int_{t}^{\infty}f_{1}\left(x\right)\,dx\\
 & \leq\int_{t}^{\infty}f_{2}\left(x\right)\,dx\\
 & =P\left(X_{2}\geq t\right).
\end{align*}
Similarly it is proved that $P\left(X_{1}\leq t\right)\geq P\left(X_{2}\leq t\right)$
for $t\leq x_{0}$ but this implies that $P\left(X_{1}>t\right)\leq P\left(X_{2}>t\right).$
If $\frac{f_{2}\left(x\right)}{f_{1}\left(x\right)}$ is increasing
then there exist a number $x_{0}$ such that $f_{1}\left(x\right)\geq f_{2}\left(x\right)$
for $x\leq x_{0}$ and that $f_{1}\left(x\right)\leq f_{2}\left(x\right)$
for $x\geq x_{0}.$

\begin{figure}
\begin{centering}
\includegraphics[scale=0.7]{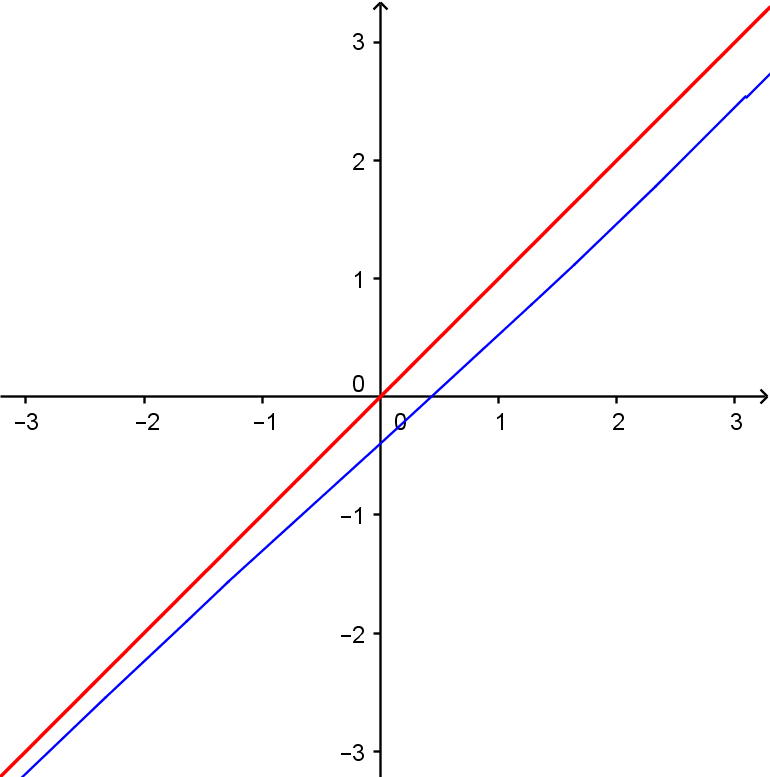}
\par\end{centering}

\caption{Plot of the quantiles of a standard Gaussian vs. the same quantiles
of the signed log-likelihood of the invers Gaussian with $\mu=1$
and $\lambda=1$ . }
\end{figure}
\end{svmultproof}

\begin{theorem}
\label{thm:InverseGaussian}If $W$ is inverse Gauiisan distributed
$IG\left(\mu,\lambda\right)$ then the signed log-likelihood 
\[
G_{W\left(\mu,\lambda\right)}\left(W\right)
\]
 is stochastically dominated by the standard Gaussian distribution,
i.e. the inequality 
\[
\Phi\left(G_{W\left(\mu,\lambda\right)}\left(w\right)\right)\leq\Pr\left(W\leq w\right)
\]
holds for any $w\in\left]0,\infty\right[$.\end{theorem}
\begin{svmultproof}
We have to prove that if $W$ has an inverse Gaussian distribution
then $G\left(W\right)$ is stochastically dominated by the standard
Gaussian. According to Lemma \ref{lem:dominans} we can prove stochastic
dominance by proving that $\phi\left(z\right)/f\left(z\right)$ is
increasing. Now
\[
\frac{\phi\left(z\right)}{f\left(z\right)}=\frac{g^{-1}\left(z\cdot\left[\frac{\mu}{\lambda}\right]^{\nicefrac{1}{2}}\right)+1}{2}
\]
which is increasing because the function $g$ is increasing. 
\end{svmultproof}

If Wald random variables are added they become more and more Gaussian
and so do their signed log-likelihood. The next theorem states that
the convergence of the signed log-likelihood towards the standard
Gaussian is monotone in stochastic domination.
\begin{theorem}
\label{Thm:IGsam}Assume that $W_{1}$ and $W_{2}$ have inverse Gaussian
distributions let $G_{1}\left(W_{1}\right)$ and $G_{2}\left(W_{2}\right)$
denote their signed log-likelihood. Then $G_{1}\left(W_{1}\right)$
is stochantically dominated by $G_{2}\left(W_{2}\right)$ if and only
if $\frac{\mu_{1}}{\lambda_{1}}>\frac{\mu_{2}}{\lambda_{2}}.$\end{theorem}
\begin{svmultproof}
We have to prove that the densities satisfy 
\[
\frac{\phi\left(z\right)2}{g^{-1}\left(z\cdot\left[\frac{\mu_{1}}{\lambda_{1}}\right]^{\nicefrac{1}{2}}\right)+1}<\frac{\phi\left(z\right)2}{g^{-1}\left(z\cdot\left[\frac{\mu_{2}}{\lambda_{2}}\right]^{\nicefrac{1}{2}}\right)+1}
\]
for $z>0$ and the reverse inequality for $z<0.$ The inequality is
equivalent to 
\[
g^{-1}\left(z\cdot\left[\frac{\mu_{1}}{\lambda_{1}}\right]^{\nicefrac{1}{2}}\right)>g^{-1}\left(z\cdot\left[\frac{\mu_{2}}{\lambda_{2}}\right]^{\nicefrac{1}{2}}\right).
\]
For $z>0$ this follows because the function $g^{-1}$ is increasing.
The reversed inequality is proved in the same way.
\end{svmultproof}

\section{Exponential distributions\label{sec:Exponential-distributions}}

Although the tail probabilities of the exponential distribution are
easy to calculate the inequalities related to the signed log-likelihood
of the exponential distribution are non-trivial and will be useful
later.

The exponential distribution $Exp^{\theta}$ has density
\[
f\left(x\right)=\frac{1}{\theta}\exp\left(-\frac{x}{\theta}\right).
\]
The distribution function is
\[
\Pr\left(X\leq x\right)=\int_{0}^{x}\frac{1}{\theta}\exp\left(-\frac{t}{\theta}\right)\,\mathrm{d}t=1-\exp\left(-\frac{x}{\theta}\right).
\]
The mean of the exponential distribution $Exp^{\theta}$ is $\theta$
and the variance is $\theta^{2}$ so the variance function is $V\left(\mu\right)=\mu^{2}.$
The divergence can be calculated as
\begin{align*}
D\left(Exp^{\theta_{1}}\Vert Exp^{\theta_{2}}\right) & =\int_{\theta_{1}}^{\theta_{2}}\frac{\mu-\theta_{1}}{\mu^{2}}\,\mathrm{d}\mu\\
 & =\frac{\theta_{1}}{\theta_{2}}-1-\ln\frac{\theta_{1}}{\theta_{2}}.
\end{align*}

From this we see that 
\begin{align*}
G_{Exp^{\theta}}\left(x\right) & =\pm\left[2\left(\frac{x}{\theta}-1-\ln\frac{x}{\theta}\right)\right]^{\nicefrac{1}{2}}\\
 & =\gamma\left(\frac{x}{\theta}\right)
\end{align*}
where $\gamma$ denotes the function 
\[
\mbox{\ensuremath{\gamma\left(x\right)}=}\begin{cases}
-\left[2\left(x-1-\ln x\right)\right]^{\nicefrac{1}{2}}, & \textrm{when}\,x\leq1;\\
+\left[2\left(x-1-\ln x\right)\right]^{\nicefrac{1}{2}}, & \textrm{when}\,x>1.
\end{cases}
\]
Note that the \emph{saddle-point approximation} is exact for the family
of exponential distributions, i.e.
\begin{align*}
f\left(x\right) & =\frac{\tau^{\nicefrac{1}{2}}}{\mathrm{e}}\cdot\frac{\phi\left(G\left(x\right)\right)}{\left[V\left(x\right)\right]^{\nicefrac{1}{2}}}.
\end{align*}

\begin{figure}
\begin{centering}
\includegraphics[scale=0.8]{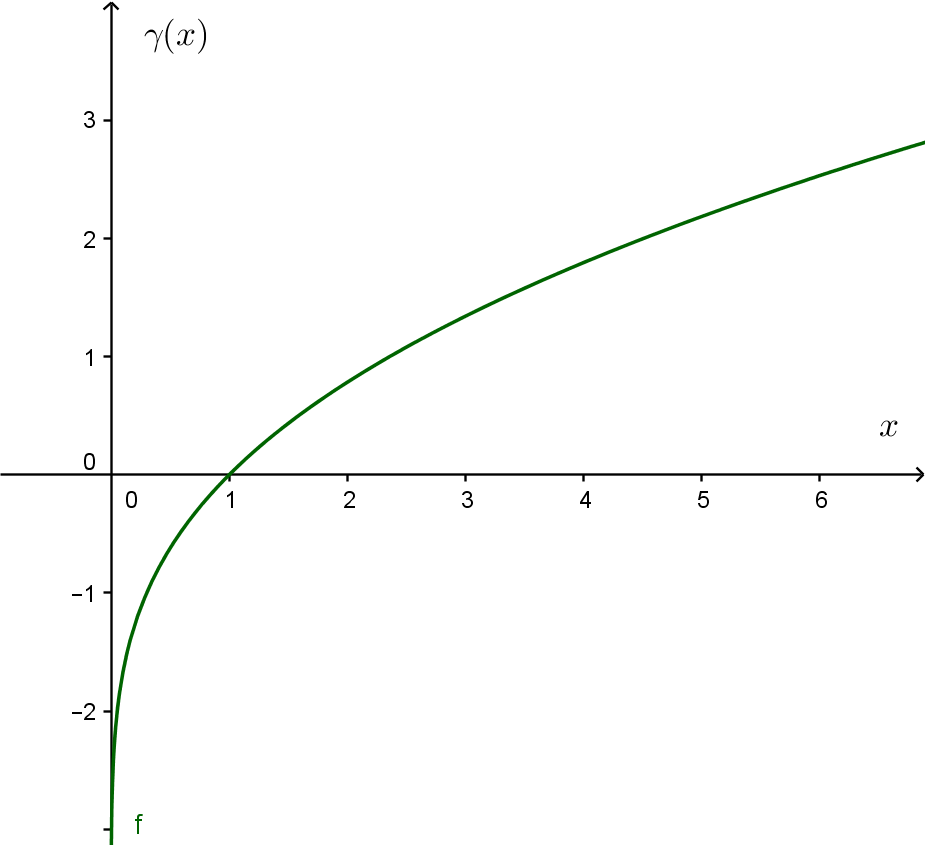}
\par\end{centering}

\caption{The signed log-likelihood $\gamma\left(x\right)$ of an exponential
distribution.}
\end{figure}

\begin{lemma}
\label{Lem:DensityExpo}The density of the signed log-likelihood of
an exponential random variable is given by 
\[
\frac{\tau^{\nicefrac{1}{2}}}{\mathrm{e}}\cdot\frac{z\phi\left(z\right)}{\gamma^{-1}\left(z\right)-1}.
\]
\end{lemma}
\begin{svmultproof}
Let $X$ be a $Exp^{\theta}$ distributed random variable. The density
of the signed log-likelihood is
\begin{align*}
\frac{f\left(G^{-1}\left(z\right)\right)}{G'\left(G^{-1}\left(z\right)\right)} & =\frac{\frac{\tau^{\nicefrac{1}{2}}}{\mathrm{e}}\cdot\frac{\phi\left(G\left(G^{-1}\left(z\right)\right)\right)}{\left[V\left(G^{-1}\left(z\right)\right)\right]^{\nicefrac{1}{2}}}}{G'\left(G^{-1}\left(z\right)\right)}\\
 & =\frac{\tau^{\nicefrac{1}{2}}}{\mathrm{e}}\cdot\frac{\phi\left(z\right)}{\left[V\left(G^{-1}\left(z\right)\right)\right]^{\nicefrac{1}{2}}G'\left(G^{-1}\left(z\right)\right)}.
\end{align*}
The variance function is $V\left(x\right)=x^{2}$ so the density is
\[
\frac{\tau^{\nicefrac{1}{2}}}{\mathrm{e}}\cdot\frac{\phi\left(z\right)}{G^{-1}\left(z\right)\cdot G'\left(G^{-1}\left(z\right)\right)}.
\]
From $G^{2}=2D$ follows that $G\cdot G'=D'$ so that 
\begin{align*}
G'\left(x\right) & =\frac{\frac{\mathrm{d}D}{\mathrm{d}x}}{G\left(x\right)}\\
 & =\frac{\frac{1}{\theta}-\frac{1}{x}}{G\left(x\right)}\,.
\end{align*}
Hence the density of $G\left(X\right)$ can be written as
\begin{align*}
\frac{\tau^{\nicefrac{1}{2}}}{\mathrm{e}}\cdot\frac{\phi\left(z\right)}{G^{-1}\left(z\right)\cdot\frac{\frac{1}{\theta}-\frac{1}{G^{-1}\left(z\right)}}{G\left(G^{-1}\left(z\right)\right)}} & =\frac{\tau^{\nicefrac{1}{2}}}{\mathrm{e}}\cdot\frac{z\phi\left(z\right)}{\frac{G^{-1}\left(z\right)}{\theta}-1}\\
 & =\frac{\tau^{\nicefrac{1}{2}}}{\mathrm{e}}\cdot\frac{z\phi\left(z\right)}{\gamma^{-1}\left(z\right)-1},
\end{align*}
which proves the lemma.\end{svmultproof}

\begin{theorem}
The signed log-likelihood of an exponentially distributed random variable
is stochastically dominated by the standard Gaussian.\end{theorem}
\begin{svmultproof}
The quotient between the density of a standard Gaussian and the density
of $G\left(X\right)$ is 
\[
\frac{\mathrm{e}}{\tau^{\nicefrac{1}{2}}}\cdot\frac{\gamma^{-1}\left(z\right)-1}{z}.
\]
We have to prove that this quotient is increasing. The function $\gamma$
is increasing so it is sufficient to prove that $\frac{t-1}{\gamma\left(t\right)}$
is increasing or equivalently that 
\[
\frac{\gamma\left(t\right)}{t-1}
\]
is decreasing. This follows from Lemma \ref{lem:Voksende} because
the variance function is increasing.
\end{svmultproof}

\begin{figure}
\begin{centering}
\includegraphics[scale=0.7]{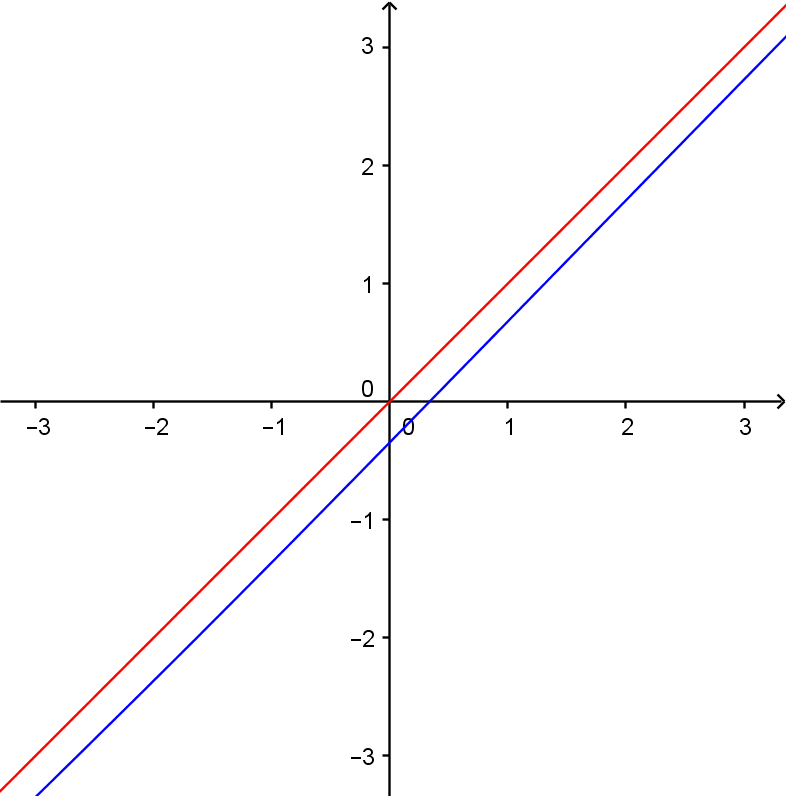}
\par\end{centering}

\caption{Plot of the quantiles of a standard Gaussian vs. the quantiles of
the signed log-likelihood of an exponential distribution.}
\end{figure}

\section{Gamma distributions\label{sec:Gamma-distributions}}

The sum of $k$ exponentially distributed random variables is Gamma
distributed $\Gamma\left(k,\theta\right)$ where $k$ is called the
shape parameter and $\theta$ is the scale parameter. It has density
\[
f\left(x\right)=\frac{1}{\theta^{k}}\frac{1}{\Gamma\left(k\right)}x^{k-1}\exp\left(-\frac{x}{\theta}\right)
\]
and this formula is used to define the Gamma distribution when $k$
is not an integer. The mean of the Gamma distribution $\Gamma\left(k,\theta\right)$
is $k\cdot\theta$ and the variance is $k\cdot\theta^{2}$ so the
variance function is $V\left(\mu\right)=\mu^{2}/k.$ The divergence
can be calculated as
\begin{align*}
D\left(\Gamma\left(k,\theta_{1}\right)\Vert\Gamma\left(k,\theta_{2}\right)\right) & =\int_{k\theta_{1}}^{k\theta_{2}}\frac{\mu-k\theta_{1}}{\mu^{2}/k}\,\mathrm{d}\mu\\
 & =k\left(\frac{\theta_{1}}{\theta_{2}}-1-\ln\frac{\theta_{1}}{\theta_{2}}\right).
\end{align*}
Further we have that 
\[
G_{\Gamma\left(k,\theta\right)}\left(x\right)=k^{\nicefrac{1}{2}}\gamma\left(\frac{x}{k\theta}\right)
\]
Note that the \emph{saddle-point approximation} is exact for the family
of Gamma distributions, i.e.
\begin{align*}
f\left(x\right) & =\frac{k^{k}\exp\left(-k\right)}{\Gamma\left(k\right)}\cdot\frac{\exp\left(-k\left(\frac{x}{k\theta}-1-\ln\frac{x}{k\theta}\right)\right)}{x}\\
 & =\frac{k^{k}\tau^{\nicefrac{1}{2}}\exp\left(-k\right)}{\Gamma\left(k\right)k^{\nicefrac{1}{2}}}\cdot\frac{\phi\left(G_{\Gamma\left(k,\theta\right)}\left(x\right)\right)}{\left[V\left(x\right)\right]^{\nicefrac{1}{2}}}.
\end{align*}

\begin{proposition}
If $F$ denotes the distribution function of the distribution $\Gamma\left(k,\theta\right)$
with mean $\mu=k\theta$ then $\frac{\mathrm{d}}{\mathrm{d}\mu}F\left(t\right)$
equals minus the density of the distribution $\Gamma\left(k+1,\theta\right)$
.\end{proposition}
\begin{svmultproof}
We have 
\begin{align*}
F\left(t\right) & =\int_{0}^{t}\frac{1}{\theta^{k}}\frac{1}{\Gamma\left(k\right)}x^{k-1}\exp\left(-\frac{x}{\theta}\right)\,\mathrm{d}x\\
 & =\int_{0}^{t/\theta}\frac{1}{\Gamma\left(k\right)}y^{k-1}\exp\left(-y\right)\,\mathrm{d}x\,.
\end{align*}
Hence
\begin{align*}
\frac{\mathrm{d}}{\mathrm{d}\theta}F\left(t\right) & =\frac{1}{\Gamma\left(k\right)}\left(\frac{t}{\theta}\right)^{k-1}\exp\left(-\frac{t}{\theta}\right)\left(-\frac{t}{\theta^{2}}\right)\\
 & =-k\frac{1}{\theta^{k+1}}\frac{1}{\Gamma\left(k+1\right)}x^{k}\exp\left(-\frac{x}{\theta}\right)\\
\frac{\mathrm{d}}{\mathrm{d}\mu}F\left(t\right) & =-\frac{1}{\theta^{k+1}}\frac{1}{\Gamma\left(k+1\right)}x^{k}\exp\left(-\frac{x}{\theta}\right)\,.
\end{align*}
The dependence on on shape and scaling is determined from the equation
\begin{align*}
D\left(\Gamma\left(k,\frac{x}{k}\right)\Vert\Gamma\left(k,\theta\right)\right) & =k\left(\frac{x}{k\theta}-1-\ln\frac{x}{k}\right)\\
 & =\frac{x}{\theta}-k-k\ln\frac{x}{k\theta}.
\end{align*}
From this we see that 
\begin{align*}
G_{k,\theta}\left(x\right) & =\pm\left[2k\left(\frac{x}{k\theta}-1-\ln\frac{x}{k\theta}\right)\right]^{\nicefrac{1}{2}}\\
 & =\pm k^{\nicefrac{1}{2}}\cdot\left[2\left(\frac{x}{k\theta}-1-\ln\frac{x}{k\theta}\right)\right]^{\nicefrac{1}{2}}\\
 & =k^{\nicefrac{1}{2}}\cdot\gamma\left(\frac{x}{k\theta}\right)
\end{align*}
which proves the proposition.
\end{svmultproof}

The following lemma is proved in the same way as Lemma \ref{Lem:DensityExpo}.
\begin{lemma}
The density of the signed log-likelihood of a Gamma random variable
is given by 
\[
\frac{k^{k}\tau^{\nicefrac{1}{2}}\exp\left(-k\right)}{\Gamma\left(k\right)k^{\nicefrac{1}{2}}}\cdot\frac{\frac{z}{k^{\nicefrac{1}{2}}}\phi\left(z\right)}{\gamma^{-1}\left(\frac{z}{k^{\nicefrac{1}{2}}}\right)-1}.
\]
\end{lemma}
\begin{theorem}
\label{thm:Gammabound}(From \cite{Harremoes2012}) The signed log-likelihood
of a Gamma distributed random variable is stochastically dominated
by the standard Gaussian, i.e. 
\[
\Pr\left(X\leq x\right)\geq\Phi\left(G_{\Gamma}\left(x\right)\right).
\]
\end{theorem}
\begin{svmultproof}
This is proved in the same way as the corresponding result for exponential
distributions.

\begin{figure}
\begin{centering}
\includegraphics[scale=0.8]{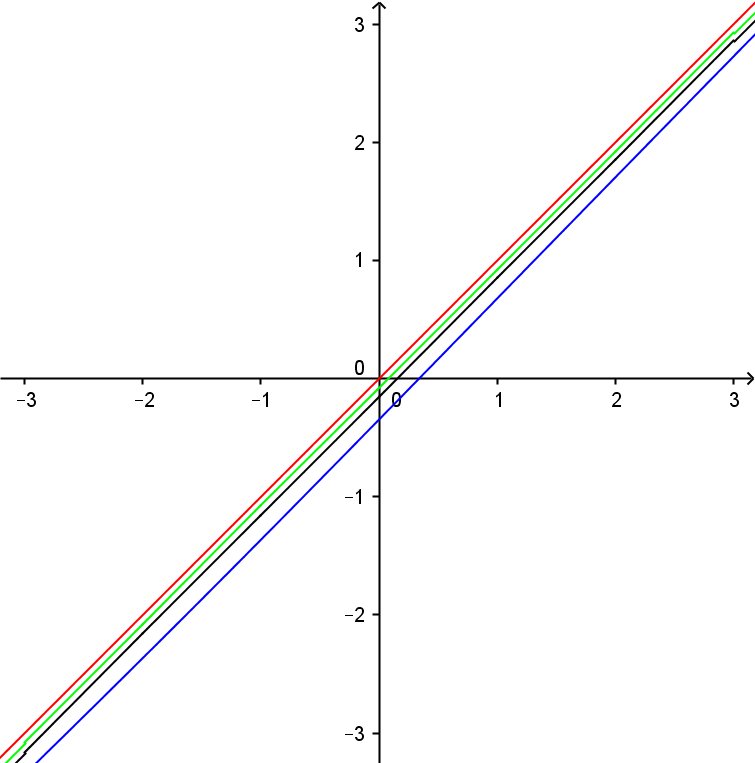}
\par\end{centering}

\caption{The quantiles of a standard Gaussian vs. gamma distributions for $k=1$
(blue), $k=5$ (black), and $k=20$ (green). The red line corresponds
to a perfect match with a standard Gaussian.}

\end{figure}
\end{svmultproof}

\begin{theorem}
\label{Thm:GammaGamma}Let $X_{1}$ and $X_{2}$ denote Gamma disstributed
random variables with shape parameters $k_{1}$ and $k_{2}$ and scale
parameters $\theta_{1}$ and $\theta_{2}.$ The the signed log-likelihood
of $X_{1}$ is dominated by the signed log-likelihood of $X_{2}$
if and only if $k_{1}\leq k_{2}.$\end{theorem}
\begin{svmultproof}
We have to prove that 
\[
\frac{\frac{z}{k_{1}^{\nicefrac{1}{2}}}\phi\left(z\right)}{\gamma^{-1}\left(\frac{z}{k_{1}^{\nicefrac{1}{2}}}\right)-1}<\frac{\frac{z}{k_{2}^{\nicefrac{1}{2}}}\phi\left(z\right)}{\gamma^{-1}\left(\frac{z}{k_{2}^{\nicefrac{1}{2}}}\right)-1}
\]
for $z>0$ and the reverse inequality for $z<0.$ The inequality is
equivalent to 
\[
\frac{\gamma^{-1}\left(\frac{z}{k_{2}^{\nicefrac{1}{2}}}\right)-1}{\frac{z}{k_{2}^{\nicefrac{1}{2}}}}<\frac{\gamma^{-1}\left(\frac{z}{k_{1}^{\nicefrac{1}{2}}}\right)-1}{\frac{z}{k_{1}^{\nicefrac{1}{2}}}}.
\]
This follows because the function 
\[
\frac{\gamma^{-1}\left(t\right)-1}{t}
\]
is increasing and because both sides have the same limit as $z$ tends
to zero from the right.
\end{svmultproof}

\section{Geometric distributions\label{sec:Geometric-distributions}}

Compounding a Poisson distribution $Po\left(\lambda\right)$ with
rate parameter $\lambda$ distributed according to an exponential
distribution $Exp\left(\theta\right)$ leads a \emph{geometric distribution}
that we will denote $Geo^{\theta}.$ We note that this is an unusual
way of parametrizing the geometric distributions, but it will be usuful
for some of our calculations. Since $\lambda$ is both the mean and
the variance of $Po\left(\lambda\right)$ the mean of $Geo^{\theta}$
is $\theta$ and the variance is $V$$\left(\mu\right)=\mu+\mu^{2}.$ 

The point probabilities of a negative binomial distribution can be
written as
\begin{align*}
\Pr\left(M=m\right) & =\int_{0}^{\infty}\frac{\lambda^{m}}{m!}\exp\left(-\lambda\right)\cdot\frac{1}{\theta}\exp\left(-\frac{\lambda}{\theta}\right)\,\mathrm{d}\lambda\\
 & =\int_{0}^{\infty}\frac{\left(\theta t\right)^{m}}{m!}\exp\left(-\theta t\right)\cdot\exp\left(-t\right)\,\mathrm{d}t\\
 & =\frac{\theta^{m}}{\left(\theta+1\right)^{m+1}}.
\end{align*}
The distribution function can be calculated as 
\begin{eqnarray*}
\Pr\left(M\leq m\right) & = & \sum_{j=0}^{m}\frac{\theta^{j}}{\left(\theta+1\right)^{j+1}}\\
 & = & 1-\left(\frac{\theta}{\theta+1}\right)^{m+1}.
\end{eqnarray*}
The divergence is given by
\begin{align*}
D\left(\left.Geo^{\theta_{1}}\right\Vert Geo^{\theta_{2}}\right) & =\int_{\theta_{1}}^{\theta_{2}}\frac{\mu-\theta_{1}}{\mu+\mu^{2}}\,\mathrm{d}\mu\\
 & =\theta_{1}\ln\frac{\theta_{1}}{\theta_{2}}-\left(\theta_{1}+1\right)\ln\frac{\theta_{1}+1}{\theta_{2}+1}.
\end{align*}
Hence the signed log-likelihood of the geometric distribution with
mean $\theta$ is given by 
\begin{align}
g_{\theta}\left(x\right) & =\pm\left[2\left(x\ln\frac{x}{\theta}-\left(x+1\right)\ln\frac{x+1}{\theta+1}\right)\right]^{\nicefrac{1}{2}}.\label{eq:gthetaDef}
\end{align}

\begin{figure}
\begin{centering}
\includegraphics[scale=0.8]{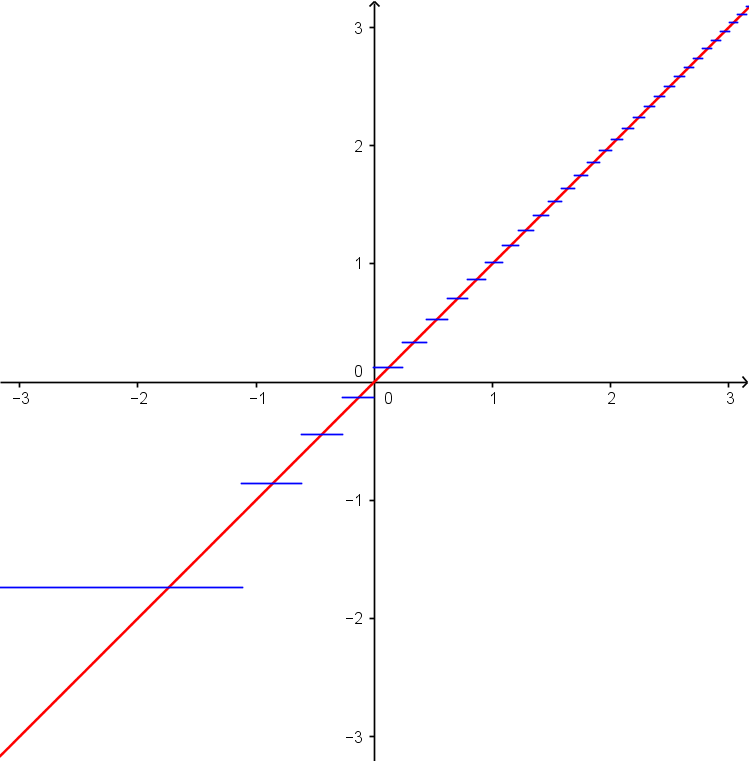}
\par\end{centering}

\caption{\label{fig:expoVsGeo}Plot the quantiles of the signed log-likelihood
of $Exp^{3.5}$ vs. the quantiles of the signed log-likelihood of
$Geo^{3.5}$.}
\end{figure}

\begin{theorem}
Assume that the random variable $M$ has a geometric distribution
$Geo^{\theta}$ and let the random variable $X$ be exponentially
distributed $Exp^{\theta}.$ If
\[
\Pr\left(X\leq x\right)=\Pr\left(M<m\right)
\]
then 
\[
G_{Geo^{\theta}}\left(m-\nicefrac{1}{2}\right)\leq G_{Exp^{\theta}}\left(x\right)\leq G_{Geo^{\theta}}\left(m\right)
\]
\end{theorem}
\begin{svmultproof}
First we note that $G_{Exp^{\theta}}\left(x\right)=\gamma\left(\nicefrac{x}{\theta}\right)$
and $\Pr\left(X\leq x\right)=\Pr\left(\nicefrac{X}{\theta}\leq\nicefrac{x}{\theta}\right).$
Therefore we introduce the variable $y=\nicefrac{x}{\theta}$ and
the random variable $Y=\nicefrac{X}{\theta}$ that is exponentially
distributed $Exp^{1}.$

We will prove that 
\begin{equation}
\Pr\left(Y\leq y\right)=\Pr\left(M<m\right)\label{eq:tailligning}
\end{equation}
 implies 
\[
g_{\theta}\left(m-\nicefrac{1}{2}\right)\leq\gamma\left(y\right)\leq g_{\theta}\left(m\right).
\]
The two inequalities are proved separately.

First we prove that $\Pr\left(Y\leq y\right)=\Pr\left(M<m\right)$
implies that $g_{\theta}\left(m-\nicefrac{1}{2}\right)\leq\gamma\left(y\right)$.
Equivalently we have to prove that 
\[
\gamma\left(y\right)-g_{\theta}\left(m-\nicefrac{1}{2}\right)=\frac{\gamma\left(y\right)^{2}-g_{\theta}\left(m-\nicefrac{1}{2}\right)^{2}}{\gamma\left(y\right)+g_{\theta}\left(m-\nicefrac{1}{2}\right)}
\]
is positive. The probability $\Pr\left(M<m\right)$ is a decreasing
function of $\theta.$ Therefore the probability $\Pr\left(Y\leq y\right)$
is a decreasing function of $\theta$, but the destribution of $Y$
does not depend on $\theta$ so $y$ must be a decreasing function
of $\theta.$ Therefore the denominator $\gamma\left(y\right)+g_{\theta}\left(m-\nicefrac{1}{2}\right)$
is a decreasing function of $\theta$ and it equals zero when $\theta=m-\nicefrac{1}{2}.$
The numerator also equals zero when $\theta=m-\nicefrac{1}{2}$ so
it is sufficient to prove that the numerator is a decreasing function
of $\theta.$ Therefore we have to prove the inequality 
\[
\frac{\partial}{\partial\theta}\left(\gamma\left(y\right)^{2}-g_{\theta}\left(m-\frac{1}{2}\right)^{2}\right)\leq0
\]
or, equivalently, that 
\[
\frac{\partial}{\partial\theta}\left(g_{\theta}\left(m-\frac{1}{2}\right)^{2}\right)\geq\frac{\partial}{\partial\theta}\left(\gamma\left(y\right)^{2}\right).
\]

One also have to prove that $\Pr\left(Y\leq y\right)=\Pr\left(M<m\right)$
implies that $\gamma\left(y\right)\leq g_{\theta}\left(m\right)$
and it is sufficient to prove that 
\[
\frac{\partial}{\partial\theta}\left(\gamma\left(y\right)^{2}\right)\geq\frac{\partial}{\partial\theta}\left(g_{\theta}\left(m\right)^{2}\right).
\]

We have 
\begin{eqnarray*}
\frac{\partial}{\partial\theta}\left(\gamma\left(y\right)^{2}\right) & = & \frac{\mathrm{d}y}{\mathrm{d}\theta}\cdot\frac{\mathrm{d}}{\mathrm{d}y}\left(\gamma\left(y\right)^{2}\right)\\
 & = & \frac{\mathrm{d}y}{\mathrm{d}\theta}\cdot2\left(1-\frac{1}{y}\right)\,.
\end{eqnarray*}
For the geometric distribution we have 
\begin{eqnarray*}
\frac{\partial}{\partial\theta}\left(g_{\theta}\left(m-1\right)^{2}\right) & = & \frac{\partial}{\partial\theta}\left(2\left(\left(m-1\right)\ln\frac{m-1}{\theta}-m\cdot\ln\frac{m}{\theta+1}\right)\right)\\
 & = & 2\left(-\frac{m-1}{\theta}+\frac{m}{\theta+1}\right)\\
 & = & 2\frac{\theta-m+1}{\theta+\theta^{2}}\,.
\end{eqnarray*}

Therefore we have to prove that
\[
2\frac{\theta-m+\nicefrac{1}{2}}{\theta+\theta^{2}}\geq2\frac{\mathrm{d}y}{\mathrm{d}\theta}\cdot\left(1-\frac{1}{y}\right)\geq2\frac{\theta-m}{\theta+\theta^{2}}\,.
\]

Therefore Equation (\ref{eq:tailligning}) can be solved as
\begin{align*}
1-\exp\left(-y\right) & =1-\left(\frac{\theta}{\theta+1}\right)^{m}\\
y & =m\ln\left(\frac{\theta+1}{\theta}\right).
\end{align*}
The derivative is 
\begin{align*}
\frac{\mathrm{d}y}{\mathrm{d}\theta} & =m\left(\frac{1}{\theta+1}-\frac{1}{\theta}\right)\\
 & =-\frac{m}{\theta+\theta^{2}}.
\end{align*}

Finally we have to prove that
\begin{eqnarray*}
\frac{\theta-m+\nicefrac{1}{2}}{\theta+\theta^{2}} & \geq & -\frac{m}{\theta+\theta^{2}}\cdot\left(1-\frac{1}{m\ln\left(\frac{\theta+1}{\theta}\right)}\right)\geq\frac{\theta-m}{\theta+\theta^{2}}\\
\theta-m+\nicefrac{1}{2} & \geq & -m+\frac{1}{\ln\left(\frac{\theta+1}{\theta}\right)}\geq-m+\theta\\
\theta+\nicefrac{1}{2} & \geq & \frac{1}{\ln\left(\frac{\theta+1}{\theta}\right)}\geq\theta\\
\left(\theta+\nicefrac{1}{2}\right)\ln\left(\frac{\theta+1}{\theta}\right) & \geq & 1\geq\theta\ln\left(1+\frac{1}{\theta}\right).
\end{eqnarray*}
The right ineqality is trivial. The left inequality is equivalent
to
\[
\ln\left(\frac{\theta+1}{\theta}\right)-\frac{1}{\theta+\nicefrac{1}{2}}\geq0.
\]
We have 
\[
\ln\left(\frac{\theta+1}{\theta}\right)-\frac{1}{\theta+\nicefrac{1}{2}}\to0
\]
 for $\theta\to\infty.$ The derivative is negative 
\[
\frac{1}{\theta+1}-\frac{1}{\theta}+\frac{1}{\left(\theta+\nicefrac{1}{2}\right)^{2}}=\frac{-1}{\theta\left(\theta+1\right)\left(2\theta+1\right)^{2}},
\]
 which proves the theorem.\end{svmultproof}

\begin{corollary}
\label{thm:ExpVsGeo}Assume that the random variable $M$ has a geometric
distribution $Geo^{\theta}$ and let the random variable $X$ be exponential
distributed $Exp^{\theta}.$ If
\[
G_{Exp^{\theta}}\left(x\right)=G_{Geo^{\theta}}\left(m\right)
\]
then
\[
\Pr\left(M<m\right)\leq\Pr\left(X\leq x\right)\leq\Pr\left(M\leq m\right).
\]

\end{corollary}
If we plot quantiles of an exponential distribution against the corresponding
quantiles of the signed log-likelihood of a Geometric distribution
we get a stair case function, i.e. a sequence of horisontal lines.
The inequality means that the left endpoint of any step is to the
left of the line $y=x.$ Actually the line $y=x$ intersects each
step and we say that the plot has an \emph{intersection property}
as illustrated in Figure \ref{fig:expoVsGeo}.
\begin{svmultproof}
Since
\[
\Pr\left(X\leq x\right)=\Pr\left(M<m\right)
\]
implies 
\[
G_{Exp^{\theta}}\left(x\right)\leq G_{Geo^{\theta}}\left(m\right)
\]
and both $\Pr\left(X\leq x\right)$ and $G_{Exp^{\theta}}\left(x\right)$
are increasing functions of $x$ we have that 
\[
G_{Exp^{\theta}}\left(x\right)=G_{Geo^{\theta}}\left(m\right)
\]
 implies that $\Pr\left(X\leq x\right)\geq\Pr\left(M<m\right).$

Since 
\[
\Pr\left(X\leq x\right)=\Pr\left(M<m\right)
\]
implies 
\[
G_{Geo^{\theta}}\left(m-\nicefrac{1}{2}\right)\leq G_{Exp^{\theta}}\left(x\right)
\]
we have that $G_{Geo^{\theta}}\left(m-\nicefrac{1}{2}\right)=G_{Exp^{\theta}}\left(x\right)$
implies that $\Pr\left(X\leq x\right)\leq\Pr\left(M<m\right)$. Hence
$G_{Geo^{\theta}}\left(m+\nicefrac{1}{2}\right)=G_{Exp^{\theta}}\left(x\right)$
implies that $\Pr\left(X\leq x\right)\leq\Pr\left(M<m+1\right)=\Pr\left(M\leq m\right).$
Since $G_{Geo^{\theta}}\left(m\right)\leq G_{Geo^{\theta}}\left(m+\nicefrac{1}{2}\right)$
we also have that $G_{Geo^{\theta}}\left(m\right)=G_{Exp^{\theta}}\left(x\right)$
implies that $\Pr\left(X\leq x\right)\leq\Pr\left(M\leq m\right).$
\end{svmultproof}

\section{Inequalities for negative binomial distributions\label{sec:Inequalities-for-waiting}}

Compounding a Poisson distribution $Po\left(\lambda\right)$ with
rate parameter $\lambda$ distributed according to a Gamma distribution
$\Gamma\left(k,\theta\right)$ leads a \emph{negative binomial distribution.}
The link to waiting times in Bernoulli processes will be explored
in Sectoin \ref{sec:Inequalities-for-binomial}. In this section we
will parametrize the negative binomial distribution as $neg\left(k,\theta\right)$
where $k$ and $\theta$ are the prameters of the corresponding Gamma
distribution. We note that this is an unusual way of parametrizing
the negative binomial distribution, but it will be usuful for some
of our calculations. Since $\lambda$ is both the mean and the variance
of $Po\left(\lambda\right)$ we can calculate the mean of $neg\left(k,\theta\right)$
as $\mu=k\theta$ and the variance as $V$$\left(\mu\right)=\mu+\frac{\mu^{2}}{k}.$ 

The point probaiblities of a negative binomial distribution can be
written in several ways
\begin{align*}
\Pr\left(M=m\right) & =\int_{0}^{\infty}\frac{\lambda^{m}}{m!}\exp\left(-\lambda\right)\cdot\frac{1}{\theta^{k}}\frac{1}{\Gamma\left(k\right)}\lambda^{k-1}\exp\left(-\frac{\lambda}{\theta}\right)\,\mathrm{d}\lambda\\
 & =\int_{0}^{\infty}\frac{\left(\theta t\right)^{m}}{m!}\exp\left(-\theta t\right)\cdot\frac{1}{\Gamma\left(k\right)}t^{k-1}\exp\left(-t\right)\,\mathrm{d}t\\
 & =\frac{\Gamma\left(m+k\right)}{m!\Gamma\left(k\right)}\cdot\frac{\theta^{m}}{\left(\theta+1\right)^{m+k}}.
\end{align*}

We need an explicite formula for the divergence that is given by
\begin{align*}
D\left(\left.neg\left(k,\theta_{1}\right)\right\Vert neg\left(k,\theta_{2}\right)\right) & =\int_{k\theta_{1}}^{k\theta_{2}}\frac{\mu-k\theta_{1}}{\mu+\frac{\mu^{2}}{k}}\,\mathrm{d}\mu\\
 & =k\left(\theta_{1}\ln\frac{\theta_{1}}{\theta_{2}}-\left(\theta_{1}+1\right)\ln\frac{\theta_{1}+1}{\theta_{2}+1}\right).
\end{align*}
The log-likelihood is given by 
\begin{align*}
G_{neg\left(k,\theta\right)}\left(x\right) & =k^{\nicefrac{1}{2}}g_{\theta}\left(\frac{x}{k}\right)
\end{align*}
where $g_{\theta}$ is given by Equation \ref{eq:gthetaDef}.

We will need the following lemma.
\begin{lemma}
A Poisson random variable $K$ with distribution $Po\left(\lambda\right)$
satisfies
\[
\frac{\mathrm{d}}{\mathrm{d}\lambda}\Pr\left(K\leq k\right)=-\Pr\left(K=k\right).
\]
\end{lemma}
\begin{svmultproof}
If $X$ is a Gamma distributed $\Gamma\left(k+1,1\right)$ then 
\begin{align*}
\Pr\left(K\leq k\right) & =\Pr\left(K<k+1\right)\\
 & =1-\Pr\left(X\leq\lambda\right).
\end{align*}
Hence
\begin{align*}
\frac{\mathrm{d}}{\mathrm{d}\lambda}\Pr\left(K\leq k\right) & =-\frac{1}{1^{k}}\frac{1}{\Gamma\left(k+1\right)}\lambda^{\left(k+1\right)-1}\exp\left(-\frac{\lambda}{1}\right)\\
 & =-\frac{\lambda^{k}}{k!}\exp\left(-\lambda\right),
\end{align*}
which proves the lemma.\end{svmultproof}

\begin{lemma}
\label{lem:difNeg}If the distribution of $M_{k}$ is $neg\left(k,\theta\right)$
then the partial derivative of the point probability with respect
to the mean value parameter equals
\[
\frac{\mathrm{d}}{\mathrm{d}\mu}\Pr\left(M_{k}\leq m\right)=-\Pr\left(M_{k+1}=m\right).
\]
where $M_{k+1}$ is $neg\left(k+1,\theta\right).$\end{lemma}
\begin{svmultproof}
We have 
\begin{align*}
\frac{\mathrm{d}}{\mathrm{d}\mu}\Pr\left(M_{k}\leq m\right) & =\frac{1}{\frac{\mathrm{d}\mu}{\mathrm{d}\theta}}\cdot\frac{\mathrm{d}}{\mathrm{d}\theta}\left(\int_{0}^{\infty}\left(\sum_{j=0}^{m}Po\left(\theta t;j\right)\right)\cdot\frac{1}{\Gamma\left(k\right)}t^{k-1}\exp\left(-t\right)\,\mathrm{d}t\right)\\
 & =\frac{1}{k}\cdot\int_{0}^{\infty}\left(-t\cdot Po\left(\theta t;m\right)\right)\cdot\frac{1}{\Gamma\left(k\right)}t^{k-1}\exp\left(-t\right)\,\mathrm{d}t\\
 & =-\int_{0}^{\infty}Po\left(\theta t;m\right)\cdot\frac{1}{\Gamma\left(k+1\right)}t^{k}\exp\left(-t\right)\,\mathrm{d}t.
\end{align*}
The last integral equals $-\Pr\left(M_{k+1}=m\right),$ which proves
the lemma.
\end{svmultproof}

\begin{flushleft}
The following theorem generalizes Corollary \ref{thm:ExpVsGeo} from
$k=1$ to arbitrary positive values of $k.$ We cannot use the same
proof technique because we do not have an explicite formula for the
quantile function for the Gamma distributions except in the case when
$k=1.$ Lemma \ref{lem:dominans} cannot be used because we want to
compare a discrete distribution with a continuous function. Instead
the proof combines a proof method developed by Zubkov and Serov \cite{Zubkov2013}
with the ideas and results developed in the previous sections.
\par\end{flushleft}
\begin{theorem}
\label{thm:NegbinVsGamma}Assume that the random variable $M$ has
a negative binomial distribution $neg\left(k,\theta\right)$ and let
the random variable $X$ be Gamma distributed $\Gamma\left(k,\theta\right).$
If
\[
G_{\Gamma\left(k,\theta\right)}\left(x\right)=G_{neg\left(k,\theta\right)}\left(m\right)
\]
then
\begin{equation}
\Pr\left(M<m\right)\leq\Pr\left(X\leq x\right)\leq\Pr\left(M\leq m\right).\label{eq:nedop}
\end{equation}
\end{theorem}
\begin{svmultproof}
Below we only give the proof of the upper bound in Inequality \ref{eq:nedop}.
The lower bound is proved the in the same way.

First we note that $G_{\Gamma\left(k,\theta\right)}\left(x\right)=G_{\Gamma\left(k,1\right)}\left(\nicefrac{x}{\theta}\right)$
and 
\[
\Pr\left(X\leq x\right)=\Pr\left(\nicefrac{X}{\theta}\leq\nicefrac{x}{\theta}\right).
\]
 Therefore we introduce the variable $y=\nicefrac{x}{\theta}$ and
the random variable $Y=\nicefrac{X}{\theta}$ that is Gamma distributed
$\Gamma\left(k,1\right).$ Introduce the difference
\[
\delta\left(\mu_{0}\right)=\Pr\left(M\leq m\right)-\Pr\left(Y\leq y\right)
\]
and note that 
\begin{equation}
\delta\left(0\right)=\lim_{\mu_{0}\to\infty}\delta\left(\mu_{0}\right)=0.\label{eq:graensevaerdier}
\end{equation}
We note that there exists (at least) one value of $\mu_{0}$ such
that $\frac{\partial\delta}{\partial\mu_{0}}=0.$ It is sufficient
to prove that $\delta$ is first increasing and then decreasing in
$\left[0,\infty\right[.$

According to Lemma \ref{lem:difNeg} the derivative of $\Pr\left(M\leq m\right)$
with respect to $\mu_{0}$ is
\begin{align*}
\frac{\partial}{\partial\mu_{0}}\Pr\left(M\leq m\right) & =-\frac{\Gamma\left(m+k+1\right)}{m!\Gamma\left(k+1\right)}\cdot\frac{\theta^{m}}{\left(\theta+1\right)^{m+k+1}}.\\
 & =-\frac{m+k}{k\left(\theta+1\right)}\frac{\Gamma\left(m+k\right)}{m!\Gamma\left(k\right)}\cdot\frac{\theta^{m}}{\left(\theta+1\right)^{m+k}}\\
 & =-\frac{m+k}{\mu_{0}+k}\Pr\left(M=m\right)\\
 & =-\frac{\hat{\theta}+1}{\theta+1}\Pr\left(M=m\right)
\end{align*}
where $\theta=\mu_{0}/k$ is the scale parameter and where and $\hat{\theta}=\nicefrac{m}{k}$
is the maximum likelihood estimate of the scale parameter. Let $\hat{\Pr}$
denote the probability of $M$ calculated with respect to this maximum
likelihood estimate $\hat{\theta}$. Then we have 
\begin{align*}
\frac{\partial}{\partial\theta}\Pr\left(M\leq m\right) & =-\frac{m+k}{\theta+1}\exp\left(-D\right)\hat{\Pr}\left(M=m\right).
\end{align*}

The condition 
\[
G_{\Gamma\left(k,\theta\right)}\left(x\right)=G_{neg\left(k,\theta\right)}\left(m\right)
\]
 can be written as
\[
k^{\nicefrac{1}{2}}\gamma\left(\frac{y}{k}\right)=k^{\nicefrac{1}{2}}g_{\theta}\left(\hat{\theta}\right)
\]
which implies
\[
\left(\gamma\left(\frac{y}{k}\right)\right)^{2}=\left(g_{\theta}\left(\hat{\theta}\right)\right)^{2}.
\]
Differentiation with respect to $\theta$ gives 
\[
2\left(1-\frac{k}{y}\right)\frac{1}{k}\frac{dy}{d\theta}=2\frac{\theta-\hat{\theta}}{\theta+\theta^{2}}
\]
so that 
\[
\frac{dy}{d\theta}=\frac{1}{\frac{1}{k}-\frac{1}{y}}\cdot\frac{\theta-\hat{\theta}}{\theta+\theta^{2}}
\]
Therefore 
\begin{align*}
\frac{\partial}{\partial\theta}\Pr\left(Y\leq y\right) & =f\left(y\right)\cdot\frac{dy}{d\theta}\\
 & =\frac{k^{k}\tau^{\nicefrac{1}{2}}\exp\left(-k\right)}{\Gamma\left(k\right)k^{\nicefrac{1}{2}}}\cdot\frac{\exp\left(-D\right)}{y}\cdot\frac{1}{\frac{1}{k}-\frac{1}{y}}\cdot\frac{\theta-\hat{\theta}}{\theta+\theta^{2}}\\
 & =\frac{k^{k}\tau^{\nicefrac{1}{2}}\exp\left(-k\right)}{\Gamma\left(k\right)k^{\nicefrac{1}{2}}}\cdot\frac{\exp\left(-D\right)}{\frac{y}{k}-1}\cdot\frac{\theta-\hat{\theta}}{\theta+\theta^{2}}\\
 & =\frac{k^{k}\tau^{\nicefrac{1}{2}}\exp\left(-k\right)}{\Gamma\left(k\right)k^{\nicefrac{1}{2}}}\cdot\frac{\exp\left(-D\right)}{\gamma^{-1}\left(g_{\theta}\left(\hat{\theta}\right)\right)-1}\cdot\frac{\theta-\hat{\theta}}{\theta\left(\theta+1\right)}
\end{align*}
Combining these results we get
\begin{multline*}
\frac{\partial\delta}{\partial\theta}=-\frac{m+k}{\theta+1}\hat{\Pr}\left(M=m\right)\cdot\exp\left(-D\right)\\
-\frac{k^{k}\tau^{\nicefrac{1}{2}}\exp\left(-k\right)}{\Gamma\left(k\right)k^{\nicefrac{1}{2}}}\cdot\frac{\exp\left(-D\right)}{\left(\gamma^{-1}\left(g_{\theta}\left(\hat{\theta}\right)\right)-1\right)}\cdot\frac{\theta-\hat{\theta}}{\theta\left(1+\theta\right)}\\
=\frac{k^{k}\tau^{\nicefrac{1}{2}}\exp\left(-k\right)}{\Gamma\left(k\right)}\cdot\frac{\exp\left(-D\right)}{\theta+1}\cdot\\
\left(\frac{\hat{\theta}-\theta}{\theta\cdot\left(\gamma^{-1}\left(g_{\theta}\left(\hat{\theta}\right)\right)-1\right)}-\frac{\Gamma\left(k\right)\left(m+k\right)}{k^{k}\tau^{\nicefrac{1}{2}}\exp\left(-k\right)}\cdot\hat{\Pr}\left(M=m\right)\right).
\end{multline*}
Remark that the first factor is positive and that 
\[
\frac{\Gamma\left(k\right)k^{\nicefrac{1}{2}}}{k^{k}\tau^{\nicefrac{1}{2}}\exp\left(-k\right)}\cdot k\left(\hat{\theta}+1\right)\cdot\hat{\Pr}\left(M=m\right)
\]
is a positive number that does not depend on $\theta.$ Therefore
it is sufficient to prove that $\frac{\hat{\theta}-\theta}{\theta\cdot\left(\gamma^{-1}\left(g_{\theta}\left(\hat{\theta}\right)\right)-1\right)}$
is a decreasing function of $\theta$, or, equivalently, to prove
that $\frac{\theta\cdot\gamma^{-1}\left(g_{\theta}\left(\hat{\theta}\right)\right)-\theta}{\hat{\theta}-\theta}$
is an increasing function of $\theta.$

The partial derivative with respect to $\theta$ is
\begin{multline*}
\frac{\left(\hat{\theta}-\theta\right)\left(\gamma^{-1}\left(g_{\theta}\left(\hat{\theta}\right)\right)+\frac{\theta\cdot\gamma^{-1}\left(g_{\theta}\left(\hat{\theta}\right)\right)}{\gamma^{-1}\left(g_{\theta}\left(\hat{\theta}\right)\right)-1}\cdot\frac{\theta-\hat{\theta}}{\theta+\theta^{2}}-1\right)+\theta\cdot\gamma^{-1}\left(g_{\theta}\left(\hat{\theta}\right)\right)-\theta}{\left(\hat{\theta}-\theta\right)^{2}}\\
=\frac{\frac{\gamma^{-1}\left(g_{\theta}\left(\hat{\theta}\right)\right)-1}{\gamma^{-1}\left(g_{\theta}\left(\hat{\theta}\right)\right)}-\frac{\left(\hat{\theta}-\theta\right)^{2}}{\hat{\theta}\left(1+\theta\right)\left(\gamma^{-1}\left(g_{\theta}\left(\hat{\theta}\right)\right)-1\right)}}{\frac{\left(\hat{\theta}-\theta\right)^{2}}{\hat{\theta}\cdot\gamma^{-1}\left(g_{\theta}\left(\hat{\theta}\right)\right)}}\,.
\end{multline*}
We have to prove that 
\begin{align*}
\frac{\gamma^{-1}\left(g_{\theta}\left(\hat{\theta}\right)\right)-1}{\gamma^{-1}\left(g_{\theta}\left(\hat{\theta}\right)\right)} & \geq\frac{\left(\hat{\theta}-\theta\right)^{2}}{\hat{\theta}\left(1+\theta\right)\left(\gamma^{-1}\left(g_{\theta}\left(\hat{\theta}\right)\right)-1\right)}\,.
\end{align*}
If $\hat{\theta}\geq\theta$ the inequality is equivalent to 
\[
\frac{\left(\gamma^{-1}\left(g_{\theta}\left(\hat{\theta}\right)\right)-1\right)^{2}}{\gamma^{-1}\left(g_{\theta}\left(\hat{\theta}\right)\right)}\geq\frac{\left(\hat{\theta}-\theta\right)^{2}}{\hat{\theta}\left(1+\theta\right)}
\]
If $\hat{\theta}<\theta$ the inequality is equivalent to 
\[
\frac{\left(\gamma^{-1}\left(g_{\theta}\left(\hat{\theta}\right)\right)-1\right)^{2}}{\gamma^{-1}\left(g_{\theta}\left(\hat{\theta}\right)\right)}\leq\frac{\left(\hat{\theta}-\theta\right)^{2}}{\hat{\theta}\left(1+\theta\right)}
\]
The equation $\frac{\left(s-1\right)^{2}}{s}=t$ can be solved with
respect to $x$, which gives the solutions $s=1+\frac{t}{2}\pm\frac{\left[t^{2}+4t\right]^{\nicefrac{1}{2}}}{2}.$
For $\hat{\theta}\geq\theta$ we get 
\begin{align*}
\gamma^{-1}\left(g_{\theta}\left(\hat{\theta}\right)\right) & \geq1+\frac{\frac{\left(\hat{\theta}-\theta\right)^{2}}{\hat{\theta}\left(1+\theta\right)}}{2}+\frac{\left[\left(\frac{\left(\hat{\theta}-\theta\right)^{2}}{\hat{\theta}\left(1+\theta\right)}\right)^{2}+4\frac{\left(\hat{\theta}-\theta\right)^{2}}{\hat{\theta}\left(1+\theta\right)}\right]^{\nicefrac{1}{2}}}{2}\\
 & =1+\left(\hat{\theta}-\theta\right)\frac{\hat{\theta}-\theta+\left[\left(\hat{\theta}+\theta\right)^{2}+4\hat{\theta}\right]^{\nicefrac{1}{2}}}{2\hat{\theta}\left(1+\theta\right)}
\end{align*}
For $\hat{\theta}<\theta$ we get 
\begin{align*}
\gamma^{-1}\left(g_{\theta}\left(\hat{\theta}\right)\right) & \geq1+\frac{\frac{\left(\hat{\theta}-\theta\right)^{2}}{\hat{\theta}\left(1+\theta\right)}}{2}-\frac{\left[\left(\frac{\left(\hat{\theta}-\theta\right)^{2}}{\hat{\theta}\left(1+\theta\right)}\right)^{2}+4\frac{\left(\hat{\theta}-\theta\right)^{2}}{\hat{\theta}\left(1+\theta\right)}\right]^{\nicefrac{1}{2}}}{2}\\
 & =1+\left(\hat{\theta}-\theta\right)\frac{\hat{\theta}-\theta+\left[\left(\hat{\theta}+\theta\right)^{2}+4\hat{\theta}\right]^{\nicefrac{1}{2}}}{2\hat{\theta}\left(1+\theta\right)}
\end{align*}
Since $\gamma$ is increasing we have to prove that 
\[
g_{\theta}\left(\hat{\theta}\right)\geq\gamma\left(1+\left(\hat{\theta}-\theta\right)\frac{\hat{\theta}-\theta+\left[\left(\hat{\theta}+\theta\right)^{2}+4\hat{\theta}\right]^{\nicefrac{1}{2}}}{2\hat{\theta}\left(1+\theta\right)}\right)
\]
or, equivalently, that 
\begin{multline*}
g_{\theta}\left(\hat{\theta}\right)-\gamma\left(1+\left(\hat{\theta}-\theta\right)\frac{\hat{\theta}-\theta+\left[\left(\hat{\theta}+\theta\right)^{2}+4\hat{\theta}\right]^{\nicefrac{1}{2}}}{2\hat{\theta}\left(1+\theta\right)}\right)\\
=\frac{\left\{ g_{\theta}\left(\hat{\theta}\right)\right\} ^{2}-\left\{ \gamma\left(1+\left(\hat{\theta}-\theta\right)\frac{\hat{\theta}-\theta+\left[\left(\hat{\theta}+\theta\right)^{2}+4\hat{\theta}\right]^{\nicefrac{1}{2}}}{2\hat{\theta}\left(1+\theta\right)}\right)\right\} ^{2}}{g_{\theta}\left(\hat{\theta}\right)+\gamma\left(1+\left(\hat{\theta}-\theta\right)\frac{\hat{\theta}-\theta+\left[\left(\hat{\theta}+\theta\right)^{2}+4\hat{\theta}\right]^{\nicefrac{1}{2}}}{2\hat{\theta}\left(1+\theta\right)}\right)}
\end{multline*}
is positive. Both the denominator and the numerator are zero when
$\theta=\hat{\theta}.$ Therefore it is sufficient to prove that both
the denominator and the numerator are decreasing functions of $\theta.$

First we prove that the denominator is decreasing. The first term
is obviously decreasing. The second term is composed of $\gamma$,
which is increasing, and $y\curvearrowright1+\frac{y}{2}\pm\frac{\left[y^{2}+4y\right]^{\nicefrac{1}{2}}}{2}$
which is increasing or decreasing depending on the sign of $\pm,$
and the function $\theta\curvearrowright\frac{\left(\hat{\theta}-\theta\right)^{2}}{\hat{\theta}\left(1+\theta\right)}$
which is decreasing when $\theta\leq\hat{\theta}$ and increasing
when $\theta\geq\hat{\theta}.$ Therefore the composed function is
a decreasing function of $\theta.$

The numerator can be written as 
\begin{multline*}
2\left\{ \hat{\theta}\ln\frac{\hat{\theta}}{\theta}-\left(\hat{\theta}+1\right)\ln\frac{\hat{\theta}+1}{\theta+1}\right\} -2\left\{ \begin{array}{c}
\left(\hat{\theta}-\theta\right)\frac{\hat{\theta}-\theta+\left[\left(\hat{\theta}+\theta\right)^{2}+4\hat{\theta}\right]^{\nicefrac{1}{2}}}{2\hat{\theta}\left(1+\theta\right)}\\
-\ln\left(1+\left(\hat{\theta}-\theta\right)\frac{\hat{\theta}-\theta+\left[\left(\hat{\theta}+\theta\right)^{2}+4\hat{\theta}\right]^{\nicefrac{1}{2}}}{2\hat{\theta}\left(1+\theta\right)}\right)
\end{array}\right\} .
\end{multline*}
We calculate the derivative with respect to $\theta$, which can be
written as 
\[
\frac{-4\frac{2\theta+\hat{\theta}+4}{\theta+\theta^{2}}\left(\theta-\hat{\theta}\right)^{2}}{\theta\left(\theta+\hat{\theta}+2\right)\left[\left(\hat{\theta}+\theta\right)^{2}+4\hat{\theta}\right]^{\nicefrac{1}{2}}+\left(\theta+2\right)\left(\left(\hat{\theta}+\theta\right)^{2}+4\hat{\theta}\right)},
\]
which is obviously less than or equal to zero. 
\end{svmultproof}

If we want to give lower bounds and upper bounds to the tail probabilities
of a negative binomial distribution the following reformulation of
Theorem \ref{thm:NegbinVsGamma} is useful.
\begin{corollary}
Assume that the random variable $M$ has a negative binomial distribution
$neg\left(k,\theta\right)$ and let the random variable $X$ be Gamma
distributed $\Gamma\left(k,\theta\right).$ If
\[
G_{\Gamma\left(k,\theta\right)}\left(x\right)=G_{neg\left(k,\theta\right)}\left(m\right)
\]
Then
\begin{equation}
\Pr\left(X\leq x_{m}\right)\leq\Pr\left(M\leq m\right)\leq\Pr\left(X\leq x_{m+1}\right)\label{eq:nedop-1}
\end{equation}
where $x_{m}$ and $x_{m+1}$ are determined by
\begin{align*}
G_{\Gamma\left(k,\theta\right)}\left(x_{m}\right) & =G_{neg\left(k,\theta\right)}\left(m\right),\\
G_{\Gamma\left(k,\theta\right)}\left(x_{m+1}\right) & =G_{neg\left(k,\theta\right)}\left(m+1\right).
\end{align*}

\end{corollary}

\section{Inequalities for binomial distributions and Poisson distributions\label{sec:Inequalities-for-binomial}}

We will prove that intersections results for binomial distributions
and Poisson distributions follows from the corresponding intersection
result for negative binomial distributions and Gamma distributions.
We note that the point probabilities of a negative binomial distribution
can be written as 
\[
\frac{\Gamma\left(m+k\right)}{m!\Gamma\left(k\right)}\cdot\frac{\theta^{m}}{\left(\theta+1\right)^{m+k}}=\frac{k^{\bar{m}}}{m!}p^{k}\left(1-p\right)^{m}
\]
 where $p=\frac{1}{1+\theta}$ and where $\bar{m}$ denotes the raising
factorial. Let $nb\left(p,k\right)$ denote a negativ binomial distribution
with succes probability $p$. Then $nb\left(p,k\right)$ is the distribution
of the number of failures before the $k$'th success in a Bernoulli
process with success probability $p.$ 

\begin{figure}
\begin{centering}
\includegraphics[scale=0.7]{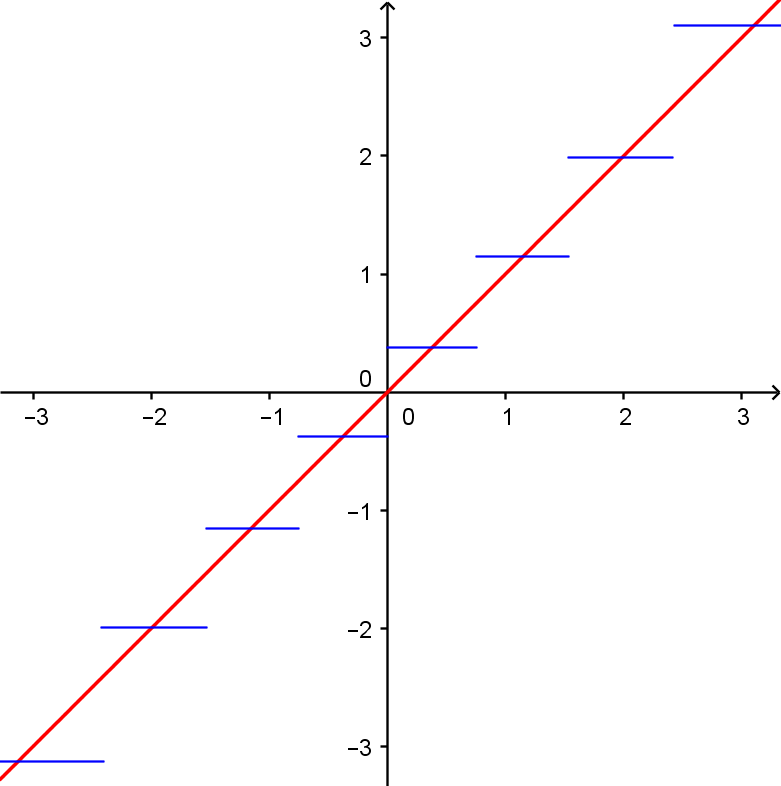}
\par\end{centering}

\caption{Plot the quantiles of the signed log-likelihood of a standard Gaussian
vs. the quantiles of the sigend log-likelihood of $bin\left(7,\nicefrac{1}{2}\right)$.}
\end{figure}

Our inequality for the negative binomial distribution can be translated
into an inequality for the binomial distribution. Assume that $K$
is binomial $bin\left(n,p\right)$ and $M$ is negative binomial $nb\left(p,k\right).$
Then
\[
\Pr\left(K\geq k\right)=\Pr\left(M+k\leq n\right).
\]
In terms of $p$ the divergence can be written as 
\begin{align*}
D\left(\left.nb\left(p_{1},k\right)\right\Vert nb\left(p_{2},k\right)\right) & =\frac{k}{p_{1}}\left(p_{1}\ln\frac{p_{1}}{p_{2}}+\left(1-p_{1}\right)\ln\frac{1-p_{1}}{1-p_{2}}\right).
\end{align*}
We have 
\[
D\left(\left.bin\left(n,p_{1}\right)\right\Vert bin\left(n,p_{2}\right)\right)=p_{1}\ln\frac{p_{1}}{p_{2}}+\left(1-p_{1}\right)\ln\frac{1-p_{1}}{1-p_{2}}
\]
so
\begin{align*}
D\left(\left.nb\left(\frac{k}{n},k\right)\right\Vert nb\left(p,k\right)\right) & =n\left(\frac{k}{n}\ln\frac{\frac{k}{n}}{p_{2}}+\left(1-\frac{k}{n}\right)\ln\frac{1-\frac{k}{n}}{1-p_{2}}\right)\\
 & =D\left(\left.bin\left(n,\frac{k}{n}\right)\right\Vert bin\left(n,p\right)\right).
\end{align*}
If $G_{bin}$ is the signed log-likelihood of $bin\left(n,p\right)$
and $G_{nb}$ is the signed log-likelihood of $nb\left(p,k\right)$
then $G_{bin\left(n,p\right)}\left(k\right)=-G_{nb\left(p,k\right)}\left(n-k\right).$ 

If $K$ is Poisson distributed with mean $\lambda$ and $X$ is Gamma
distributed with shape parameter $k$ and scale parameter 1, i.e.
the distribution of the waiting time until $k$ observations from
an Poisson process with intensity 1. Then
\[
\Pr\left(K\geq k\right)=\Pr\left(X\leq\lambda\right).
\]

Next we note that 
\[
D\left(Po\left(k\right)\Vert Po\left(\lambda\right)\right)=D\left(\left.\Gamma\left(k,\frac{\lambda}{k}\right)\right\Vert \Gamma\left(k,1\right)\right).
\]
If $G_{Po\left(\lambda\right)}$ is the signed log-likelihood for
$Po\left(\lambda\right)$ and $G_{\Gamma\left(k,1\right)}$ is the
signed log-likelihood for $\Gamma\left(k,1\right)$ then $G_{Po\left(\lambda\right)}\left(k\right)=-G_{\Gamma\left(k,1\right)}\left(\lambda\right).$ 
\begin{theorem}
\label{thm:BinomialPoisson}Assume that $K$ is binomially distributed
$bin\left(n,p\right)$ and let $G_{bin\left(n,p\right)}$ denote the
signed log-likelihood function of the exponential family based on
$bin\left(n,p\right).$ Assume that $L$ is a Poisson random variable
with distribution $Po\left(\lambda\right)$ and let $G_{Po\left(\lambda\right)}$
denote the signed log-likelihood function of the exponential family
based on $Po\left(\lambda\right).$ If 
\[
G_{bin\left(n,p\right)}\left(k\right)=G_{Po\left(\lambda\right)}\left(k\right)
\]
Then
\begin{equation}
\Pr\left(K<k\right)\leq\Pr\left(L<k\right)\leq\Pr\left(K\leq k\right).\label{eq:BinoPoisson}
\end{equation}
\end{theorem}
\begin{svmultproof}
Let $M$ denote a negative binomial random variable with distribution
$nb\left(p,k\right)$ and let $X$ denote a Gamma random variable
with distribution $\Gamma\left(k,\theta\right)$ where the parameter
$\theta$ equals $\frac{1}{p}-1$ such that the distributions $nb\left(p,k\right)$
and $\Gamma\left(k,\theta\right)$ have the same mean value. Now $G_{nb\left(p,k\right)}\left(n-k\right)=-G_{bin\left(n,p\right)}\left(k\right)$
and $G_{\Gamma\left(k,\theta\right)}\left(\lambda\theta\right)=-G_{Po\left(\lambda\right)}\left(k\right).$
Then $G_{nb\left(p,k\right)}\left(n-k\right)=G_{\Gamma\left(k,\theta\right)}\left(\lambda\theta\right).$
The left part of the Inequality \ref{eq:BinoPoisson} is proved as
follows.
\begin{align*}
\Pr\left(K<k\right) & =1-\Pr\left(K\geq k\right)\\
 & =1-\Pr\left(M+k\leq n\right)\\
 & \leq1-\Pr\left(X\leq\lambda\theta\right)\\
 & =1-\Pr\left(L\geq k\right)\\
 & =\Pr\left(L<k\right).
\end{align*}
The right part of the inequality is proved in the same way.
\end{svmultproof}

Note that Theorem \ref{thm:NegbinVsGamma} cannot be proved from Theorem
\ref{thm:BinomialPoisson} because the number parameter for a binomial
distribution has to be an integer while the number parameter of a
negative binomial distribution may assume any positive value. Now,
our inequalities for negative binomial distributions can be translated
into inequalities for binomial distributions.

\begin{figure}
\begin{centering}
\includegraphics[scale=0.7]{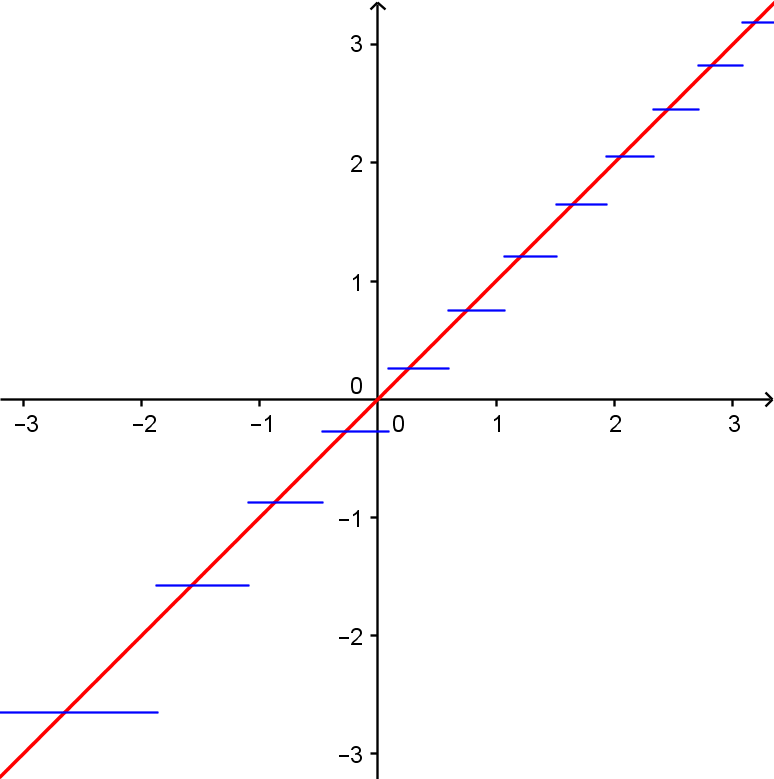}
\par\end{centering}

\caption{\label{fig:QQGausPoisson}Plot of quantiles of a standard Gaussian
vs. the log-ligelihood of the Poisson distribution $Po\left(3.5\right).$}
\end{figure}

Now we can prove the an \emph{intersection inequalities} for the binomial
family as stated in the following theorem that was recently proved
by Serov and Zubkov \cite{Zubkov2013}.
\begin{corollary}
\label{thm:Binomialintersection}Assume that $K$ is binomially distributed
$bin\left(n,p\right)$ and let $G_{bin\left(n,p\right)}$ denote the
signed log-likelihood function of the exponential family based on
$bin\left(n,p\right).$ Then
\begin{equation}
\Pr\left(K<k\right)\leq\Phi\left(G_{bin\left(n,p\right)}\left(k\right)\right)\leq\Pr\left(K\leq k\right).\label{eq:binoInter}
\end{equation}

Similarly, assume that $L$ is Poisson distributed $Po\left(\lambda\right)$
and let $G_{Po\left(\lambda\right)}$ denote the signed log-likelihood
function of the exponential family based on $Po\left(\lambda\right).$
Then
\begin{equation}
\Pr\left(L<k\right)\leq\Phi\left(G_{Po\left(\lambda\right)}\left(k\right)\right)\leq\Pr\left(L\leq k\right).\label{eq:PoissonInter}
\end{equation}
\end{corollary}
\begin{svmultproof}
First we prove the left part of Inequality (\ref{eq:PoissonInter}).
Let $X$ denote a Gamma distributed $\Gamma\left(k,1\right)$ and
let $Z$ denote a standard Gaussian. Then $G_{Po\left(\lambda\right)}\left(k\right)=-G_{\Gamma\left(k,1\right)}\left(\lambda\right)$
and 
\begin{align*}
\Pr\left(L<k\right) & =1-\Pr\left(L\geq k\right)\\
 & =1-\Pr\left(X\leq\lambda\right)\\
 & =\Pr\left(X\geq\lambda\right)\\
 & \leq\Pr\left(Z\geq G_{\Gamma\left(k,1\right)}\left(\lambda\right)\right)\\
 & =\Pr\left(Z\geq-G_{Po\left(\lambda\right)}\left(k\right)\right)\\
 & =\Phi\left(G_{Po\left(\lambda\right)}\left(k\right)\right).
\end{align*}

The left part of Inequality (\ref{eq:binoInter}) is obtained by combining
the left part of Inequality (\ref{eq:PoissonInter}) with the left
part of Inequality (\ref{eq:BinoPoisson}).\end{svmultproof}

\begin{proof}
The right part of Inequality (\ref{eq:binoInter}) is obtained follows
from the left part of Inequality (\ref{eq:binoInter}) by replacing
$p$ by $1-p$ and replacing $k$ by $n-k.$ \end{proof}

\begin{svmultproof}
Since a Poisson distribution is a limit of binomial distributions
the right part of Inequality (\ref{eq:PoissonInter}) follows from
the right part of Inequality (\ref{eq:PoissonInter}).
\end{svmultproof}

The intersection property for Poisson distributions was proved in
\cite{Harremoes2012} where the inequality for binomial distributions
was also conjectured. 

\pagebreak{}

\section{Summary}

The main theorems in this paper are domination theorems and intersection
theorems. The first type of inequalities states that the signed log-likelihood
of one distribution is dominated by the signed loglikelihood of another
distribution, i.e. the distribution function of the first distribution
is larger than the distribution function of the second distribution.

\begin{table}[h]
\begin{centering}
\begin{tabular}{|c|c|c|c|}
\hline 
signed ll & dom. by signed ll & Condition & Theorem\tabularnewline
\hline 
\hline 
Inverse Gaussian & Gaussian &  & \ref{thm:InverseGaussian}\tabularnewline
\hline 
$IG\left(\mu_{1},\lambda_{1}\right)$ & $IG\left(\mu_{2},\lambda_{2}\right)$ & $\frac{\mu_{1}}{\lambda_{1}}>\frac{\mu_{2}}{\lambda_{2}}$ & \ref{Thm:IGsam}\tabularnewline
\hline 
Gamma & Gaussian &  & \ref{thm:Gammabound}\tabularnewline
\hline 
$\Gamma_{k_{1},\theta_{1}}$ & $\Gamma_{k_{2},\theta_{2}}$ & $k_{1}\leq k_{2}$ & \ref{Thm:GammaGamma}\tabularnewline
\hline 
\end{tabular}
\par\end{centering}

\caption{Stochastic domination results. Note that the exponential distributions
are special cases of Gamma distributions.}
\end{table}
The second type of result are intersection results, i.e. the distribution
function of the log-likelihood of a discrete distribution is a staircase
function where each step is intersected by the distribution function
of the log-likelihood of a continuous distribution. 

\begin{table}[h]
\begin{centering}
\begin{tabular}{|c|c|c|}
\hline 
Discrete distribution & Continuous distribution & Theorem\tabularnewline
\hline 
\hline 
Geometric & Exponential & \ref{thm:ExpVsGeo}\tabularnewline
\hline 
Negative binomial & Gamma & \ref{thm:NegbinVsGamma}\tabularnewline
\hline 
Binomial & Gaussian & \ref{thm:Binomialintersection}\tabularnewline
\hline 
Poisson & Gaussian & \ref{thm:Binomialintersection}\tabularnewline
\hline 
\end{tabular}
\par\end{centering}

\caption{Intersection results.}
\end{table}

\section{Discussion\label{sec:Discussion}}

\begin{figure}
\begin{centering}
\includegraphics[scale=0.8]{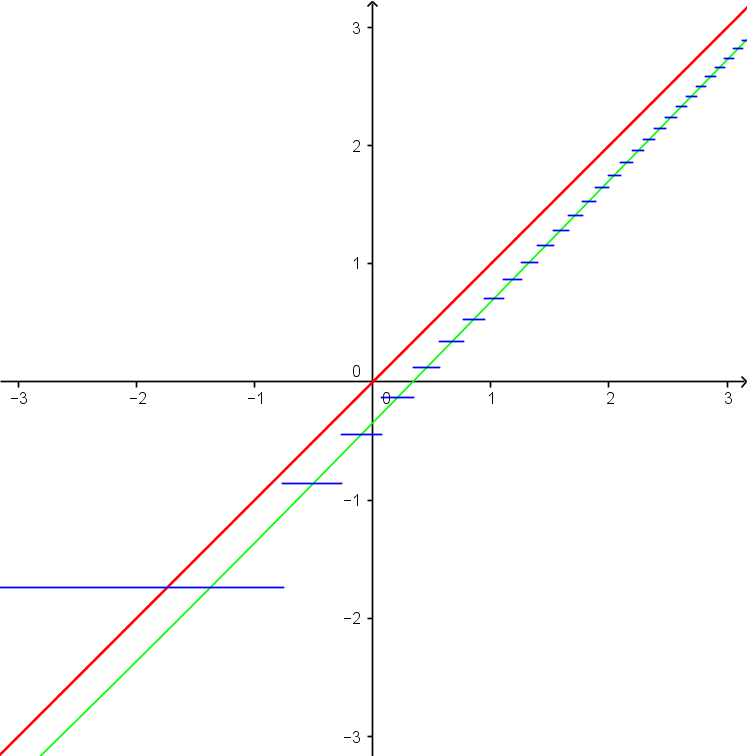}
\par\end{centering}

\caption{\label{fig:NegGammaGauss}Plot of the quantiles of a standard Gaussian
vs. the quantiles of the signed log-likelihood of the negative binomial
distribution $neg\left(1,3.5\right)$ (blue) and of the signed log-likelihood
of the Gamma distribution $\Gamma\left(1,3.5\right)$ (green).}
\end{figure}

In this paper we have presented inequalities of two types. The inequalities
for inverse Gaussian distributions, exponential distributions and
other Gamma distributions are about stochastic domination. The inequalities
for Poisson distributions, binomial distributions, geometric distributions
and other negative binomial distributions are about intersection.
These inequalities can be combined in order to get inequalities of
other types. For instance a negative binomial random variable $M$
with distribution $neg\left(k,\theta\right)$ satisfies
\[
\Phi\left(G_{neg\left(k,\theta\right)}\left(m\right)\right)\leq\Pr\left(M\leq m\right),
\]
where $G_{nb\left(p,k\right)}$ denotes the signed log-likelihood
of the negative binomial distribution. Contrary to the similar inequality
for the binomial distribution the inquality $\Pr\left(M<m\right)\leq\Phi\left(G_{neg\left(k,\theta\right)}\left(m\right)\right)$
does in general \emph{not} hold as illustrated in Figure \ref{fig:NegGammaGauss}. 

We have both lower bounds and upper bounds on the Poission distributions.
The upper bound for the Poisson distribution corresponds to the lower
bound for the Gamma distribution presented in Theorem \ref{thm:Gammabound},
but the lower for the Poisson distribution translated into a new upper
bound for the distribution function of the Gamma distribution. Numerical
calculations also indicates that in Inequality (\ref{eq:PoissonInter})
the right hand inequality can be improved to
\[
\Phi\left(G_{Po\left(\lambda\right)}\left(k+\frac{1}{2}\right)\right)\leq\Pr\left(L\leq k\right).
\]
This inequality is much tighter than the inequality in (\ref{eq:PoissonInter}).
Similarly, J. Reiczigel, L. Rejt\H{o} and G. Tusn{\'a}dy conjectured that
both the lower bound and the upper bound in Inequality \ref{eq:binoInter}
can be significantly improved when for $p=\nicefrac{1}{2}$ \cite{Reiczigel2011},
and their conjecture has been a major motivation for initializing
this research.

For the most important distributions like the binomial distributions,
the Poisson distributions, the negative binomial distributions, the
inverse Gaussian distributions and the Gamma distributions we can
formulate sharp inequalities that hold for any sample size. All these
distributions have variance functions that are polynomials of order
2 and 3. Natural exponential families with polynomial variance functions
of order at most 3 have been classified \cite{Morris1982,Letac1990}
and there is a chance that one can formulate and prove sharp inequalities
for each of these exponential families. Although there may exist very
nice results for the rest of the exponential families with simple
variance functions the rest of these exponential families have much
fewer applications than the exponential families that have been the
subject of the present paper.

In the present paper inequalities have been developed for specific
exponential families, but one may hope that some more general inequalities
may be developed where bounds on the tails are derived directly from
the properties of the variace function. 
\begin{acknowledgements}
The author want to thank Gabor Tusn{\'a}dy, who inspired me to look at
this type of inequalities. The author also want the thank L{\'a}szl{\'o} Gy{\"o}rfi,
Joram Gat, Jan{\'a}s Komlos, and A. Zubkov for useful correspondence or
discussions. Finally I want to express my gratityde to Narayana Prasad
Santhanam who invitet me to a two month visit at the Electrical Engineering
Department, University of Hawai'i. This paper was completed during
my visit.
\end{acknowledgements}


\end{document}